\documentclass[11pt]{article}

\usepackage{amsthm}
\usepackage{authblk}
\usepackage[T1]{fontenc}        
\usepackage[a4paper]{geometry}
\usepackage[latin1]{inputenc}   
\usepackage{vmargin}            
\usepackage{fancyhdr}      
\usepackage{calc}               
\usepackage[dvips]{graphicx}
\usepackage{epsf}
\usepackage[dvips]{epsfig}
\usepackage[dvipsnames]{xcolor}
\usepackage{enumitem}
\usepackage{pifont} 
\usepackage{fancybox}           
\usepackage[Glenn]{fncychap}
\usepackage{amsmath}    
\usepackage{amssymb}    
\usepackage{amsfonts}   
\usepackage{dsfont}
\usepackage{bclogo}
\usepackage[most]{tcolorbox} 
\usepackage{mathabx}

\usepackage[colorlinks=true]{hyperref}
\hypersetup{urlcolor=[rgb]{0.00,0.38,0.00},
linkcolor=[rgb]{1.00,0.00,0.20},
citecolor=[rgb]{0.85,0.00,0.00},
colorlinks=TRUE
}

\newcommand{\lsup}[1]{\overline{#1}}
\newcommand{\linf}[1]{\underline{#1}}
\newcommand{\diame}[2]{|{#1}|_{\mathsf{#2}}}
\newcommand{\conf}[3]{{\vphantom{{#2}_{j}}}_{\footnotesize{{{^{#1}}}}}\!{#2}_{#3}}
\def\segN#1#2{{[\![#1,#2]\!]}} 

\newcommand{\pto}[1]{\ensuremath{\text{o}\left ( #1 \right )}}
\newcommand{\gdo}[1]{\ensuremath{O\left ( #1 \right )}}
\DeclareMathOperator*{\E}{\mathbb{E}}
\newcommand{\Chi}{\mathcal{X}}
\newcommand{\fonc}{\mathfrak{F}}
\newcommand{\mass}{\mathfrak{m}}

\newcommand{\ds}{\displaystyle}

\newcommand{\plus}{\mbox{\raisebox{.2mm}{\footnotesize{\ensuremath{+}}}}}

\newcommand{\pinf}{\plus\ensuremath{\infty}}

\newcommand{\Prod}{\ensuremath{\prod\limits}}
\newcommand{\PI}[2]{\ensuremath{%
        \mathop{\mbox{\Large\ensuremath{\Pi}}}\limits}_{#1}^{#2}}

\renewcommand{\emptyset}{\ensuremath{\varnothing}}

\newcommand{\N}{\ensuremath{\mathbb{N}}}

\newcommand{\R}{\ensuremath{\mathbb{R}}}

\newcommand{\C}{\ensuremath{\mathbb{C}}}
\newcommand{\X}{\ensuremath{\mathbb{X}}}
\newcommand{\D}{\ensuremath{\mathbb{D}}}

\newcommand{\K}{\ensuremath{\mathbb{K}}}

\renewcommand{\P}{\ensuremath{\mathcal{P}}}

\DeclareMathOperator{\diam}{diam}
\DeclareMathOperator{\dist}{{\rm dist}}

\renewcommand{\le}{\ensuremath{\leqslant}}
\renewcommand{\ge}{\ensuremath{\geqslant}}

\newcommand{\scaln}[3]{\ensuremath{\langle #1,#2 \rangle_{#3}}}

\theoremstyle{plain}
\newtheorem{Th}{Theorem}[section]
\newtheorem{Lemme}[Th]{Lemma}
\newtheorem{Cor}[Th]{Corollary} 
\newtheorem{Prop}[Th]{Proposition}
\newtheorem{Rem}[Th]{Remark}
\newtheorem{Def}[Th]{Definition}
\newtheorem*{Nota}{Notation}
\newtheorem*{Propri}{Properties}

\newcommand*{\email}[1]{%
    \normalsize\href{mailto:#1}{#1}\par
    }

\begin{document}

\title{Dimensions of harmonic measures on non-autonomous Cantor sets}
\author{Athanasios Batakis}
\affil{Institut Denis Poisson, Universit\'e d'Orl\'eans, BP 6759, 45067 Orl\'eans cedex 2, France\\
\email{athanasios.batakis@univ-orleans.fr}}
\author{Guillaume Havard}
\affil{Universit\'e d'Orl\'eans, BP 6759, 45067 Orl\'eans cedex 2, France\\
\email{guillaume.havard@univ-orleans.fr}}

\maketitle 

\begin{abstract}
We consider Non Autonomous Conformal Iterative Function Systems (NACIFS) and their limit set. Our main concern is harmonic measure and its dimensions : Hausdorff and Packing. We prove that this two dimensions are continuous under perturbations and that they verify Bowen's and Manning's type formulas. In order to do so we prove general results about measures, and more generally about positive functionals, defined on a symbolic space, developing tools from thermodynamical formalism in a non-autonomous setting.

\end{abstract}

\section{Introduction and  statement of results.}\label{definition}
Given a compact set in the plane defined dynamically by iteration of conformal contractions, this article focuses on the variations of Hausdorff dimension and packing dimension of the harmonic measure under perturbations of the conformal contractions. In a natural framework, which is precisely stated below, we demonstrate that these functions are continuous.

Let $\K\subset \C$ be a compact set and $\omega$ the harmonic measure of $\hat\C\setminus \K$  evaluated at $\infty$. Since the works of Makarov \cite{Ma}, Jones, and Wolff \cite{JW}, it is  known that the Hausdorff dimension of $\omega$ is less than or equal to $1$. This result was generalized by Bourgain to $\mathbb{R}^d$ in \cite{Bourgain}, where he proved the existence of a constant $\ell(d)>0$ such that the dimension of the harmonic measure of any open set in $\mathbb{R}^d$ is less than $d-\ell(d)$. According to Makarov, Jones, and Wolff we know that $\ell(2)=1$.

If $\K$ has Hausdorff dimension in $]1,2]$ the previous results tells us that harmonic measure only "sees" a small part of $\K$ of dimension less or equal to $1$. Let for instance $\K$ be the Mandelbrot set. This is a compact and connected set, so that $\C\setminus \K$, is simply connected and we know by \cite{Ma} that the support of harmonic measure is a subset of the boundary of $\K$ of Hausdorff dimension equal to $1$. On the other side Shishikura proved in \cite{Shishi} that any intersection of this boundary with an open set has Hausdorff dimension $2$.

When the compact set $\K$ is fractal and has Hausdorff dimension less than one it is still reasonable to believe that in many cases the dimension of harmonic measure remains stricly less than the Hausdorff dimension of $\K$. This has been proved to be true for several types of conformal Cantor sets \cite{Ca}, \cite{LyB}, \cite{UZ}, \cite{Batakis}. More recently Batakis and Zdunik proved in \cite{bazd} that this remains true for a class of non autonomous Cantor sets generated by similitudes. Note that in our knowledge the question whether this is true for all Conformal cantor sets is still opened. On the other side the first author exhibits in \cite{Batakis} examples of Cantor sets for which harmonic measure has full dimension, and very recently David, Jeznach and Julia in \cite{DJJ} give an example of a Cantor set in the plane for which  Harmonic measure is equivalent with Hausdorff measure.

In this paper, we develop tools to study entropies, Lyapunov exponents, pressures and dimensions for non-autonomous Cantor sets. Most results are proved in a symbolic space and for the set of positive functions defined on cylinders, see section \ref{symbolspace}, even if this is interesting on his own, our main focus in this article are Hausdorff and Packing dimension of the harmonic measure, and their continuity properties, for a class of non-autonomous conformal repellers in the plane. 

In order to state our results precisely, we describe now our settings and we introduce some notations.

Following \cite{ReUr}, given a Jordan domain $U$, a non-autonomous conformal iterated function system is a sequence $\Psi=({\Psi_n})$ of finite collections of conformal contractions $\Psi_{n}=\{\conf{n}{\psi}{i}:U\to U\;;\; i=1,...,d_n\}$ verifying properties (1), (2) and (3) below for some $\eta>0$
There exists $\eta>0$ such that
\begin{enumerate}
\item \underline{\em Conformality:}   There is a fixed neighborhood $V\supset(1+\eta)U$ of $U$ such that  $\conf{n}{\psi}{i}$ extends to a conformal diffeomorphism from $V$ to $V$, for all $1\le i\le d_n$ and all $n\in\N$ .

\item  \underline{\em Open set condition and annulus condition:} 
\begin{itemize} 
\item $\conf{n}{V}{i}\cap \conf{n}{V}{j}=\emptyset$ for all $1\le i\not=j\le d_n$, where $\conf{n}{V}{i}=\conf{n}{\psi}{i}(V)$. [OSC]
\item\label{AC} Moreover, for any $n$ and any $1\le i\le d_n$ we have $\conf{n}{V}{i}\subset (1-\eta)U$. [AC]
\end{itemize}
\item \underline{\em Bounded Contraction:} For all $n\in\N$, all  $1\le i\not=j\le d_n$ and all $x\in\conf{n}{V}{i}$
$$\eta < \|\conf{n}{\psi}{i}'(x) \|< 1-\eta\quad\mbox{[BC]}$$

\end{enumerate}

Combining the open set and the bounded contraction conditions we get a uniform bound on $d_n$:
\begin{equation}\label{CDCS}
\mbox{there exists } N \mbox{ s.t. } d_n\le N \mbox{  for all } n\in\N.
\end{equation}

\begin{Nota}\label{admissible_map} 
\begin{itemize}
    \item 
    For any $E\subset V$, any $n\in \N$ and any $1\le i\le d_n$ let $\conf{n}{E}{i}=\conf{n}{\psi}{i}(E)$. The index on the left stands for the {\bf generation} of the map while the one on the right stands for the chosen \textbf{branch} of the map.
    \item 
    For $n\in\N$ we denote by ${\mathcal U}_n$ the collection $\{\conf{n}{U}{1},\cdots,\conf{n}{U}{d_{n}}\}$. 
    Let  $\Phi_{n}=\ds\sum_{i=1}^{d_n}\conf{n}{\phi}{i}\mathds{1}_{{}_{\conf{n}{U}{i}}}$  where $\conf{n}{\phi}{i}$, defined on $\conf{n}{U}{i}$, is the inverse conformal dilatation of $\conf{n}{\psi}{i}$, ie. $\conf{n}{\phi}{i}:\conf{n}{U}{i}\to U$ and $\conf{n}{\phi}{i}\circ \conf{n}{\psi}{i}=\mbox{id}_{U}$.
\end{itemize}
\end{Nota}
In this way, we have introduced  piecewise conformal maps denoted again by $\Phi_n$ defined on the union of topological discs:  $\ds\Phi_n: \bigcup_{\conf{n}{U}{i}\in\, \mathcal{U}_n} \conf{n}{U}{i}\to U$ by the formula
$$\Phi_{n}(z)=\conf{n}{\phi}{i}(z)\,\,\mbox{ for}\,\, z\in \conf{n}{U}{i}.$$
\begin{Rem}
Note that the condition of bounded contraction is stronger than the one of  \cite{ReUr} (uniform contraction). This is due to the fact that dealing with harmonic measure implies estimates on capacity and hence on the size of the sub-components. Also the conformality condition has been strengthened  to ensure existence of conformal annuli with bounded from below moduli. This is crucial in order to be able to use Harnack's principle (cf. \cite{Ca}).

Furthermore, without loss of generality we can suppose that $U$ is the unit disk.
\end{Rem}
We consider the following non-autonomous limit set

\begin{Def}\label{admcantor}  Given a non-autonomous iterated conformal functions system $\Psi_n$  on a Jordan domain $U$ we define the {\em limit set}  $\X\subset \mathbb{C}$ by

$$\X=\bigcap_{n=1}^\infty \left (\Phi_n\circ \Phi_{n-1}\circ\dots \circ \Phi_2\circ \Phi_1\right )^{-1}(U),$$

where $\Phi_n$ is the sequence of piecewise conformal dilatations introduced above. 
\end{Def}
For any $n$, the map $T_{n}:=\Phi_n\circ \Phi_{n-1}\circ\dots \circ \Phi_2\circ \Phi_1$ is well defined on $\X$ and for any $x\in\X$ we have $T_{n}(x)\in U$. Moreover there exists $\PI{i=1}{n}d_{i}$ conformal inverse branches of $T_{n}$ defined on $V$. Each of this branches sends $U$ onto a conformal disk included in $U$,  each point in $\X$ belongs to one and only one of these conformal disks. Remember that  $\mathcal{U}_{n}$ is the collection of connected components of $T_n^{-1}(U)$.
\\There is an obvious way to code the situation. 
\\For any integer $n$, let $\mathcal{A}_{n}:=\{1\cdots d_{n}\}$ and set $\K:=\ds\PI{k=1}{\pinf}\mathcal{A}_{k}$. 
There is a one-to-one correspondence from $\K$ to $\X$. More precisely if $\mathtt{a}=(a_{n})$ is an element of $\K$, then the sequence of subsets of $V$ defined by $U_{0}=\bar{U}$ and $U_{n}=\conf{n}{\psi}{a_{n}}(U_{n-1})$ is a decreasing sequence of compact sets with diameters converging to $0$, thus converging to a point $x$, which is obviously an element of $\X$.
Reciprocally, if $x\in\X$ then for any $n$ we have : $\ds \Phi_{n}\circ\cdots\circ \Phi_{1}(x)\in U$, and there exists $\mathtt{a}=(a_{n})\in\K$ such that $x\in \conf{1}{U}{a_{1}}$, $\ds\Phi_{1}(x)\in\conf{2}{U}{a_{2}},\cdots, \Phi_{n}\circ\cdots\circ \Phi_{1}(x)\in\conf{{n+1}}{U}{a_{n+1}}$.
\\If $x\in\X$ is associated with $\mathtt{a}\in \K$, then for any $n$ we have $x\in \conf{1}{\psi}{_{a_{1}}}\circ\cdots\circ\conf{n}{\psi}{_{a_{n}}}(U)$, which is the conformal disk mentioned above. In the sequel we will denote $X_{n}(x):=a_{1}\cdots a_{n}$ the set $\X \cap \conf{1}{\psi}{_{a_{1}}}\circ\cdots\circ\conf{n}{\psi}{_{a_{n}}}(U)$, and $\Psi_{a_{1}\cdots a_{n}}^{-1}:=\conf{1}{\psi}{_{a_{1}}}\circ\cdots\circ\conf{n}{\psi}{_{a_{n}}}$.

Recall the definition of Hausdorff Dimension and Packing Dimension of a probability measure $\mu$:
$$HD(\mu)=\inf_{\{Z:\mu(Z)=1\}}HD(Z)\quad\mbox{and}\quad PD(\mu)=\inf_{\{Z:\mu(Z)=1\}}PD(Z),$$
the infimum and the supremum being taken over all Borel subsets $Z$ of the support of $\mu$. 
Let $\omega$ be  the harmonic measure on $\hat{\mathbb{C}}\setminus \X$ evaluated at $\infty$. By the celebrated  results of N. Makarov \cite{Ma} and P. Jones-T. Wolff \cite{JW} the Hausdorff dimension of $\omega$ is not larger than one. On the other hand, it is clear that  the Hausdorff dimension of $\omega$ is at most $HD(\X)$.

For any continuous function defined on $U$ let  $\ds\|f\|:=\sup_{U}|f|$, with $|\,\cdot\,|:\C\to \R^+$ is the modulus.
Let $\X$ and $\tilde \X$ be two non-autonomous limit sets associated to the conformal iterated systems $\Psi $ and $\tilde \Psi$ respectively. We define 
$$d(\Psi,\tilde\Psi)=\sup_n\max_{1\le i\le d_{n}}\left(\left\| \conf{n}{\psi}{i}-\conf{n}{\tilde{\psi}}{i}\right\|  \right)$$

Using general results on the symbolic space $\K_{1}$, and studied in the second and third part of this paper, we prove 
\begin{Th}\label{Continuity}
Le $\X$ and $(\X_{k})_{k\in\N}$ be non-autonomous limit sets  associated with the conformal iterated systems $\Psi $ and $(\Psi_{k})_{k\in\N}$ respectively. Let $\omega$ and $\omega_{k}$ be the harmonic measures of their complementaries. Assume that $\ds\lim_{k\to \pinf}d(\Psi,\Psi_{k})=0$, then, 
$$\lim_{k\to\pinf}HD(\omega_{k})=HD(\omega)\quad\mbox{and}\quad \lim_{k\to\pinf}PD(\omega_{k})=PD(\omega)\cdot$$
\end{Th}

We also discuss validity of well known formulas to calculate these dimensions. In particular we show that an adapted form of Manning's formula applies.

\begin{Th}\label{thdim}
Let $\omega$ be the harmonic measure of  $\C\setminus\X$, where $\X$ is the limit set of a  non-autonomous conformal iterated systems $\Psi $.

We have 
$$HD(\omega)=\liminf_{n}\frac{H_{\omega,n}}{\chi_{\omega,n}}\quad\mbox{and}\quad PD(\omega)=\limsup_{n}\frac{H_{\omega,n}}{\chi_{\omega,n}},$$
where $H_{\omega,n}=\ds-\frac1n\sum_{X\in {\mathcal U}_n}\omega(X)\log\omega(X)$ and  $\chi_{\omega,n}=\ds\frac1n\sum_{X\in{\mathcal U}_n}\omega(X)\log\|T_{n}'(X)\|$.

Moreover, $\omega$-almost surely
$$HD(\omega)=\liminf_{n}\frac{h_{\omega,n}(x)}{\chi_{n}(x)}\quad\mbox{and}\quad PD(\omega)=\limsup_{n}\frac{h_{\omega,n}(x)}{\chi_{n}(x)}$$
 where for any $x\in\X_{0}$, $h_{\omega,n}(x)=-\ds\frac{1}{n}\log \omega(X_{n}(x))$ and $\chi_{n}(x)=\ds\frac{1}{n}\log |T_{n}'(x)|\cdot$
\end{Th}

\section{Symbolic point of view}\label{symbolspace}
This paragraph deals exclusively with symbolic Cantor sets. More precisely, for any integer $n$, let $\mathcal{A}_{n}$ be a finite alphabet with $d_{n}$ symbols. Set $\K_{n}:=\ds\prod_{k=n}^{\pinf}\mathcal{A}_{k}$. We are mainly interested in $\K_{1}$ that we denote $\K$. 

 Let $T:\K_{n}\to\K_{n+1}$ be the shift map sending a sequence $\ds\mathtt{a}=(a_{i})_{i\ge 0}$ to $T(\mathtt{a})=(a_{i})_{i\ge 1}$. It should be emphasized that $T$ goes from $\K_{n}$ to $\K_{n+1}$. This is a bit ambiguous, since the same letter designates different maps, but it lightens considerably the notations. We let $\ds D^p_{n}:=\prod_{k=0}^{n-1}d_{k+p}$ be the degree of $T^n$ on $\K_{p}$.

Given a sequence of symbols $\mathtt{a}=\ds(a_{i})_{i\in\N}\in \K$,  we set 
\[X_{n}(\mathtt{a}):=\ds a_{0}\cdots a_{n-1}:=\{\mathtt{b}=(b_{i})\in \K\,|\, a_{i}=b_{i}\;, \mbox{ for } 0\le i\le n-1\}.\] 
The set $X_{n}(\mathtt{a})$ is the \textit{cylinder} (word), of \textit{length} (size) $n$,  in which $\mathtt{a}$ belongs to. Given a cylinder $X$ we denote by $|X|$ its length and, for any $n$, by $\mathcal{C}_{n}$ the set of all cylinders of length $n\ge 1$, in which, for practical reasons, we add $\emptyset$ with the convention $\emptyset X=X\emptyset=X$, and $|\emptyset|=0$. Note that  $\mathcal{C}_{n}$ is a finite partition of $\K$ by $D_n^1$ cylinders. 

The set $\K$ endowed with the product topology is compact. Let $\mathsf{d}$ be any metric on $\K$ that generates this topology and such that ${\rm diam}\K=1$. If $E\subset \K$, we denote by $|E|_{\mathsf{d}}:=\sup\{\mathsf{d}(x,y)\,|\, (x,y)\in E\times E\}$ its diameter.

Actually, since the planar Cantor sets that we will consider are conformal repellers, they are well described using the coding and the symbolic point of view. It will then be natural to choose a metric which ensures that a cylinder in $\K$, and its ``image'' in the plane, have approximately the same diameters.

\subsection{Asymptotic Siblings Invariance}
Because of the strong intra-scalar-similarity of the Cantor sets we are dealing with, one expects that natural geometric measures, or natural geometric objects, defined on them have nice scaling or invariance properties.

Note that in this section we do not need to assume any assumption on the finite numbers $d_n$.

In order to define precisely the properties we have in mind we need to introduce some definitions.

For any integer $p\ge 1$, let $\fonc_{p}$ be the set of real positive maps defined on the set of cylinders $T^p\mathcal{(C)}:=\ds\bigcup_{n\ge p}T^p\mathcal{C}_{n}$. For $p=0$ we denote by $\fonc=\fonc_{0}$ the set of maps defined on the set of cylinders $\mathcal{C}:=\ds\bigcup_{n\ge 0}\mathcal{C}_{n}$ with values in $]0,\pinf[$. 

For $\nu=\nu_{\emptyset}\in\fonc$ and $XY$ a fixed element of $\mathcal{C}$ we note $\nu_{X}(Y):=\ds\frac{\nu(XY)}{\nu(X)}$. Note that $\nu_{X}$ is an element of $\fonc_{|X|}$.

A finite measure $\nu$ on $\K_p$ that gives positive mass to any open set, is naturally an element of $\fonc_p$. Reciprocally
we say that $\nu\in \fonc_p$ is a (finite) measure if it projects as such on the compact set defined asymptotically by the cylinders, which is if and only if :
$$\forall n\in \N\quad\forall X \in T^p(\mathcal{C})\quad \frac{1}{\nu(X)}\sum_{|Y|=n}\nu(XY)=1\cdot$$

\begin{Def}\label{mass}$ \ $
\begin{itemize}
\item
Given $\nu$ and $\nu'$ in $\fonc_{p}$ we naturally define $\nu+\nu'$, $\lambda\nu$ for any $\lambda\in\R$, and $\nu.\nu'$.
\item
For any real function $\varphi$ and any $\nu\in\fonc_{p}$ we note $\varphi(\nu)$ the map from $T^p\mathcal{C}$ to $\R$ defined by $\varphi(\nu)(X)=\varphi(\nu(X))$. Clearly this only makes sense if $\varphi$ is well-defined on $\nu (T^p\mathcal{C})$.
\item
Given $\nu$ and $\nu'$ in $\fonc_{p}$ and $n\ge 1$ let : $\ds\scaln{\nu}{\nu'}{n}:=\sum_{X\in T^p\mathcal{C}_{n+p}}\nu(X)\nu'(X)\cdot$
\item
For any $\nu\in\fonc_{p}$, we define its \textbf{mass} on $\mathcal{C}_{n}$ (with $n\ge p$)
$$\mass_{n}(\nu):=\scaln{\nu}{1}{n}\cdot$$
Note that $\nu$ is a measure if $\mass_{n}(\nu_{X})=1$ for any $n$ and any $X\in\mathcal{C}$.\\
\item 
We say that $\nu$ and $\nu'$ in $\fonc_{p}$ are equivalent if there exists $C>1$ such that for any $X$ we have $\ds \frac{1}{C}\le \frac{\nu(X)}{\nu'(X)}\le C.$
We then note $\nu\sim\nu'$.
\end{itemize}
\end{Def}

In this paper we will mainly consider the following elements $\nu\in\fonc_{p}$ : 
\begin{itemize}
\item
$\nu(X)=\mu(X)$, where $\mu$ is a measure on $\K$,
\item
$\nu(X)=|X|_{\mathsf{d}}$, 
\end{itemize}

Here are some elementary properties that will be used later on. 
\begin{Propri} $\ $\\
For any $s\in \fonc$, any  $XYZ\in \mathcal{C}_{n}$ we have 
\begin{equation}
s_{X}(YZ)=s_{X}(Y)s_{XY}(Z)\qquad\mbox{and}\qquad s_{X_{Y}}(Z)=s_{XY}(Z)
\end{equation}
For any finite measure $\nu$, any cylinder $X\in\mathcal{C}$ and any integers $n$ and $p$ we have
\begin{equation}\label{equalityg}
\scaln{\nu_{X}}{\log s_{X}}{n+p}=\scaln{\nu_{X}}{\log s_{X}}{n}+\sum_{Y\in T^{|X|}\mathcal{C}_{|X|+n}}\nu_{X}(Y)\scaln{\nu_{XY}}{\log s_{XY}}{p}\cdot
\end{equation}
\end{Propri}

The following is an application of Jensen inequality.
\begin{Prop}\label{jensen}
Let $\nu$ and $\nu '$ be two elements of $\fonc_{p}$.  For any $n$ we have
\begin{equation}\label{jens1}
\ds\scaln{\nu}{\log \frac{\nu'}{\nu}}{n}\le \mathsf{m}_{n}(\nu) \log \frac{\mathsf{m}_{n}(\nu')}{\mathsf{m}_{n}(\nu)},
\end{equation}
with equality if and only if there exists $k$ such that $\nu(X)=k\nu'(X)$ for any $X\in T^p\mathcal{C}_{n+p}$.
\end{Prop}
We introduce now our main assumption we will refer to as : Asymptotic Sibling Invariance. It is a kind of long range Markov property of a functional that asserts that along cylinders which share a big common part of history, the relative functional does not depend on its origin.
\begin{Def}\label{translation}
Let $\nu$ be an element of $\fonc$ and $(\beta_{n})$ a positive sequence converging to $0$. We say that $\nu$ satisfies Asymptotic Siblings Invariance (ASI) with sequence $(\beta_{n})$, if for any integers $n$, $k$ and $p$ we have
\begin{equation}\label{ASI}
 \forall (X,X')\in \mathcal{C}_{n}^2\,\,\, \forall Y\in T^n\mathcal{C}_{n+k}\,\,\,\forall Z\in T^{n+k}\mathcal{C}_{n+k+p}\quad\left|\log\left(\frac{\nu_{XY}(Z)}{\nu_{X'Y}(Z)}\right)\right|\le \beta_{k}
\end{equation}
\end{Def}
Note that we sometimes use the following version of the inequation :
$$
\left|\frac{\nu_{XY}(Z)}{\nu_{X'Y}(Z)}-1\right|\le \beta_{k}\cdot
$$
\begin{Nota}
To simplify notations with cylinders, we sometimes are a little bit imprecise and instead of writing $XY$ with $X\in \mathcal{C}_{n}$ and $Y\in T^n\mathcal{C}_{n+k}$,
we will write $XY$ with $|X|=n$ and $|Y|=k$, the context reminding us that $Y$ is not in $\mathcal{C}_{k}$ but in $T^n\mathcal{C}_{n+k}$.
\end{Nota}

\begin{Def}\label{compatible}
Given two sequences of positif  real numbers $(\gamma_{n})$ and $(\beta_{n})$, we say that they are compatible if $(\gamma_{n})$ is increasing, $(\beta_{n})$ converges towards $0$ and there exists a sequence of integers, $(p_{n})$, such that $p_{n}\le n$ and $\ds\left(\frac{1}{n}\max (\gamma_{p_{n}}, \gamma_{n}\beta_{p_{n}})\right)$ converges towards $0$.
We then denote $\ds c(\gamma_{n},\beta_{n}):=\max (\gamma_{p_{n}}, \gamma_{n}\beta_{p_{n}})$.
\end{Def}

As an example, one may consider the situation we will have to deal with studying the harmonic measure. In that case we will have $\gamma_{n}=\gdo{n}$, it appears that such a sequence is compatible with any sequence $(\beta_{n})$ converging towards $0$. As a matter of fact  we have $\ds\frac{1}{n}\max (\gamma_{p_{n}}, \gamma_{n}\beta_{p_{n}})=\max (\frac{\gdo{p_n}}{n}, \frac{\gdo{n}}{n}\beta_{p_{n}}) $, and any  $p_{n}=\pto{n}$, for instance $p_n=\lfloor\sqrt{n}\rfloor$, fulfilled the desired growth condition.

\begin{Prop}\label{sameasymptoticvalues}
Let $(\gamma_{n})$ and $(\beta_{n})$ be two compatible sequences.
\\Assume that $s$ and $\nu$ are two elements of $\fonc$, that $\nu$ is a finite measure and that they both satisfied (ASI) with $(\beta_{n})$. 
Assume moreover that one of this condition is fulfilled :
\begin{align} \label{massexpbis1}
   &\forall XY\in\mathcal{C}\quad \left|\log s_{X}(Y)\right|\le \gamma_{|Y|}&\tag{\theequation.1}\\
&\mbox{ or }& \notag \\
\label{massexpbis2}
&\forall X\in\mathcal{C}\quad\forall n\in\N\quad s(X)\ge1\quad\mbox{and}\quad \scaln{\nu_{X}}{\log(s_{X})}{n}\le \gamma_{n}\cdot& \tag{\theequation.2}
\end{align}

Then for any $\Tilde{s}\sim s$
\begin{itemize}
\item[--]
There exists $C>0$ such that : 
$$\alpha_{n}(\nu,\Tilde{s}):=\sup_k\left\{\max_{|X|=k}\scaln{\nu_{X}}{\log \Tilde{s}_{X}}{n}-\min_{|X|=k}\scaln{\nu_{X}}{\log \Tilde{s}_{X}}{n}\right\}\le Cc(\gamma_{n},\beta_{n})\cdot$$
\item[--]
For any $(X,X')\in\mathcal{C}^2$ we have 
$$\lim_{n\to\pinf}\left(\frac{1}{n}\scaln{\nu_{X'}}{\log \Tilde{s}_{X'}}{n}- \frac{1}{n}\scaln{\nu_{X}}{\log \Tilde{s}_{X}}{n}\right)=0\cdot$$
\end{itemize}
\end{Prop}

\proof
The constant $C$ will change from line to line but will stay independent of integers or cylinders, only depending on $\nu$ and $s$.
\\
For any $\Tilde{s}\sim s$ we have, for any $XY\in \mathcal{C}$ : 
$-2\log  C \le\log \Tilde{s}_X(Y)-\log s_X(Y)\le 2\log C$.

From where we deduce, because $\nu$ is a finite measure
$$-2\log C\le\scaln{\nu_{X}}{\log \Tilde{s}_{X}}{n}-\scaln{\nu_{X}}{\log s_{X}}{n}\le 2\log C\cdot$$
This tells us that conclusions of the theorem are true for $\Tilde{s}$ as soon as they are true for $s$, and in the sequel we focus on $s$.
\\Let  
\begin{itemize}
\item[--] for any $X\in\mathcal{C}$, $g_{n}(X):=\scaln{\nu_{X}}{\log s_{X}}{n},$ with $g_{n}:=g_{n}({\emptyset})$
\item[--] $\ds\alpha_{n}:=\alpha_{n}(\nu,s)=\sup_k\left\{\max_{|X|=k}g_{n}(X)-\min_{|X|=k}g_{n}(X)\right\}\cdot$
\end{itemize}

We first note that (\ref{massexpbis1}), and the fact that $\nu$ is a finite measure, imply that for any $X\in\mathcal{C}$ we have
\begin{equation}\label{massexp}
|g_{n}(X)|=\left|\scaln{\nu_{X}}{\log s_{X}}{n}\right|\le \gamma_{n},
\end{equation}
inequality which is also true with (\ref{massexpbis2}) and leads to 
\begin{equation}\label{alphgamm}
\alpha_{n}\le 2\gamma_{n}\cdot
\end{equation}

Note also that given any $X\in\mathcal{C}$, any integers $n$ and $p$, we have from (\ref{equalityg}) :
\begin{eqnarray}\label{massexp4}
g_{n+p}(X)=g_n(X)+\sum_{|Y|=n}\nu_{X}(Y)g_{p}(XY)\cdot
\end{eqnarray}

From where we deduce that for any $X ' \in\mathcal{C}_{|X|}$ we have 

\begin{eqnarray}\label{massexp1}
g_{n+p}(X)-g_{n+p}(X') & = & g_{n}(X)-g_{n}(X')+\sum_{|Y|=n}\left(\nu_{X}(Y)-\nu_{X'}(Y)\right)g_{p}(XY)+\\\nonumber
 & & \sum_{|Y|=n}\nu_{X'}(Y)\left(g_{p}(XY)-g_{p}(X'Y)\right)
\end{eqnarray}

Let $\mathcal{B^+}$ be the subset of $T^{|X|}\mathcal{C}_{|X|+n}$ for which we have $\nu_{X}(Y)-\nu_{X'}(Y)\ge 0$, and let $\mathcal{B^-}$ be its complementary. Then we have
$$\sum_{|Y|=n}\left(\nu_{X}(Y)-\nu_{X'}(Y)\right)g_{p}(XY)\le \max_{|Z|=|X|+n}g_{p}(Z)\sum_{Y\in \mathcal{B^+}}\left(\nu_{X}(Y)-\nu_{X'}(Y)\right)+$$
$$\qquad\qquad\qquad\qquad\qquad\qquad\qquad\qquad\qquad\qquad\qquad \min_{|Z|=|X|+n}g_{p}(Z)\sum_{Y\in \mathcal{B^-}}\left(\nu_{X}(Y)-\nu_{X'}(Y)\right)\cdot$$
Since $\nu$ is a finite measure, we have $\ds \sum_{|Y|=n}\nu_{X}(Y)=\sum_{|Y|=n}\nu_{X'}(Y)=1$.\\
This implies that  $\ds \sum_{|Y|=n}\left(\nu_{X}(Y)-\nu_{X'}(Y)\right)=0$.\\
So that $\ds  \sum_{Y\in \mathcal{B^+}}\left(\nu_{X}(Y)-\nu_{X'}(Y)\right)=-\sum_{Y\in \mathcal{B^-}}\left(\nu_{X}(Y)-\nu_{X'}(Y)\right)$, which leads to
\begin{eqnarray*}
\sum_{|Y|=n}\left(\nu_{X}(Y)-\nu_{X'}(Y)\right)g_{p}(XY) & \le & \sum_{Y\in \mathcal{B^+}}\left(\nu_{X}(Y)-\nu_{X'}(Y)\right)\left(\max_{|Z|=|X|+n}g_{p}(Z)-\min_{|Z|=|X|+n}g_{p}(Z)\right)\\
 & \le &\alpha_{p} \sum_{Y\in \mathcal{B^+}}\left(\nu_{X}(Y)-\nu_{X'}(Y)\right)
 \\
 & \le &\alpha_{p} \sum_{Y\in \mathcal{B^+}}\nu_{X'}(Y)\left(\frac{\nu_{X}(Y)}{\nu_{X'}(Y)}-1\right)
 \\
 & \le &\alpha_{p}\beta_0 \sum_{Y\in \mathcal{B^+}}\nu_{X'}(Y)\qquad\mbox{(because of (ASI))}
 \\
 & \le & \alpha_{p}\beta_0 \qquad\mbox{(because $\nu$ is a finite measure)}
 \end{eqnarray*}
We thus have 
\begin{eqnarray}\label{massexp2}
\sum_{|Y|=n}\left(\nu_{X}(Y)-\nu_{X'}(Y)\right)g_{p}(XY)\le \lambda\alpha_{p}\cdot
\end{eqnarray}
\\From (ASI) we deduce that  : 
$$\ds\nu_{XY}(Z)\le \left(1+\beta_{|Y|}\right)\nu_{X'Y}(Z),$$
and also : 
$$\ds \log s_{XY}(Z)\le\log s_{X'Y}(Z) +\beta_{|Y|}\cdot$$
\\This implies :
\begin{itemize}
\item[--] if (\ref{massexpbis1}) is fulfilled, that we have $\ds |\log s_{X'Y}(Z)|\le \gamma_{|Z|}$ so that
$$\begin{array}{rcl}
\nu_{XY}(Z)\log  s_{XY}(Z)&-&\nu_{X'Y}(Z)\log s_{X'Y}(Z)\\
&=&  (\nu_{XY}(Z)-\nu_{X'Y}(Z))\log s_{XY}(Z)+\nu_{X'Y}(Z)\log \frac{s_{XY}(Z)}{s_{X'Y}(Z)}\\
& \le & |\nu_{XY}(Z)-\nu_{X'Y}(Z)| \gamma_{|Z|}+\nu_{X'Y}(Z)\left|\log \frac{s_{XY}(Z)}{s_{X'Y}(Z)}\right|\\
\end{array}$$
and since $\nu$ and $s$ fulfill (ASI) we get 
$$\begin{array}{rcl}
\phantom{\nu_{XY}(Z)\log  s_{XY}(Z)}
& \le & \beta_{|Y|}\nu_{X'Y}(Z)\gamma_{|Z|}+\nu_{X'Y}(Z)\beta_{|Y|}\qquad\qquad\qquad\qquad\quad\quad\\
& \le & \beta_{|Y|}(1+\gamma_{|Z|})\nu_{X'Y}(Z).
\end{array}$$
Summing over all cylinders $Z$ of length $p$, and using once again the fact that $\nu$ is a finite measure, we get 
\begin{eqnarray}\label{massexp6}
g_{p}(XY)-g_{p}(X'Y)\le (1+\gamma_{p})\beta_{|Y|}\le C\gamma_{p}\beta_{|Y|}\cdot
\end{eqnarray}
\item[--] if (\ref{massexpbis2}) is fulfilled, that we have
$$
\begin{array}{rcl}
\nu_{XY}(Z)\log  s_{XY}(Z)-\nu_{X'Y}(Z)\log s_{X'Y}(Z) & \le & \beta_{|Y|}\left(\nu_{X'Y}(Z)\log s_{X'Y}(Z)+C\nu_{X'Y}(Z)\right)
\end{array}
$$
Summing over all cylinders $Z$ of length $p$, and using (\ref{massexp}), we also get (\ref{massexp6}).

\end{itemize}

Injecting (\ref{massexp6}) and (\ref{massexp2}), into  (\ref{massexp1}) leads to :
$$g_{n+p}(X)-g_{n+p}(X') \le g_{n}(X)-g_{n}(X')+C\gamma_{p}\beta_{n}+\lambda\alpha_{p}\cdot$$
From where we get that : 
$$\ds \alpha_{n+p}\le \alpha_{n}+C\gamma_{p}\beta_{n}+\lambda\alpha_{p}\cdot$$
Which might be rewritten, for any $p\le n$, in the following form
$$\ds \alpha_{n}\le \alpha_{p}+C\gamma_{n-p}\beta_{p}+\lambda\alpha_{n-p}\le \alpha_{p}+C\gamma_{n}\beta_{p}+\lambda\alpha_{n-p}\,,\quad\mbox{$(\gamma_{n})$ being increasing.} $$
Using (\ref{alphgamm}) we get 
$$\ds \alpha_{n}\le 2\gamma_{p}+C\gamma_{n}\beta_{p}+\lambda \alpha_{n-p}\le C\max(\gamma_{p},\gamma_{n}\beta_{p})+\lambda \alpha_{n-p}\cdot$$
 
Let $q$ and $r\le p$ be the quotient and the remainder of the euclidean division of $n$ by $p$. By induction we get
$$\begin{array}{rcl}
\ds \alpha_{n} & \le &  \ds C\max(\gamma_{p},\gamma_{n}\beta_{p})\sum_{k=0}^{q-1} \lambda^k+\lambda^q\alpha_{r}\\
&\le &\ds C\max(\gamma_{p},\gamma_{n}\beta_{p})+2\lambda^q\gamma_{r}\\
 & \le & \ds C\max(\gamma_{p},\gamma_{n}\beta_{p})+2\lambda^q\gamma_{p}\\
 & \le & \ds C\max(\gamma_{p},\gamma_{n}\beta_{p})\cdot
\end{array}
$$
Since $(\gamma_{n})$ and $(\beta_{n})$ are compatible, one may take $p=p_{n}$ to get
$$ \alpha_{n}\le Cc(\gamma_{n},\beta_{n}),$$
which is the first part of the proposition and tells us that for any  $X$ and $X'$ in $\mathcal{C}$, with $|X|=|X'|$, we have
$$\frac{1}{n}\left|g_{n(X)}-g_{n(X')}\right|\le \frac{\alpha_{n}}{n}\le C\frac{c(\gamma_{n},\beta_{n})}{n}\to_{n\to\pinf}0\cdot$$

We show now that we can get rid of the condition $|X|=|X'|$. 

Let $X\in \mathcal{C}$, and let $a$ be a cylinder of length 1. Apply (\ref{massexp4}) with $n=1$ and $p=q$
$$g_{q+1}(X)=g_{1}(X)+\sum_{|b|=1}\nu_{X}(b)g_{q}(Xb)=g_{q}(Xa)+g_{1}(X)+\sum_{|b|=1}\nu_{X}(b)(g_{q}(Xb)-g_{q}(Xa))\cdot$$
Apply now (\ref{massexp4})  with $n=q$ and $p=1$ :
$$g_{q+1}(X)=g_{q}(X)+\sum_{|Y|=q}\nu_{X}(Y)g_{1}(XY).$$
We thus have :
$$g_{q}(X)-g_{q}(Xa)=g_{1}(X)+\sum_{|b|=1}\nu_{X}(b)(g_{q}(Xb)-g_{q}(Xa))-\sum_{|Y|=q}\nu_{X}(Y)g_{1}(XY).$$
By hypothesis on $\nu$ and $s$, there exists $C>0$ such that for any cylinder $Z\in\mathcal{C}$ we have $g_{1}(Z)\le C$ . This leads to
$$\left| g_{q}(X)-g_{q}(Xa)\right|\le C+\alpha_{q},$$
and this easily leads, for any $Y$ with $|Y|=p$, to
\begin{eqnarray}\label{ineq2}
\left| g_{q}(X)-g_{q}(XY)\right|\le (C+\alpha_{q})p\le C\alpha_{q}p\cdot
\end{eqnarray}
Let $X'\in \mathcal{C}_{|X|}$ and $XY\in\mathcal{C}$, we have
$$\left| g_{q}(X')-g_{q}(XY)\right|\le \left| g_{q}(X')-g_{q}(X)\right|+\left| g_{q}(X)-g_{q}(XY)\right|\le \alpha_{q}+Cp\alpha_{q}\le C'p\alpha_{q},$$
and then
$$\frac{1}{q}\left| g_{q}(X')-g_{q}(XY)\right|\le C'p\frac{\alpha_{q}}{q}\cdot$$
Which finishes the proof.
\qed

\subsection{A continuity result}
Let $\nu$ and $\nu'$ be two elements of $\fonc$, we define 
$$\mathcal{D}(\nu,\nu')=\sup_{p\ge0}\max\left\{\left|\log\frac{\nu_{X}(a)}{\nu'_{X}(a)}\right|,\,\,Xa\in \mathcal{C}_{p+1}\mbox{ and }|a|=1\right\}\cdot$$

Note that we obviously have $\mathcal{D}(\nu,\nu')=\mathcal{D}(\nu',\nu)$ and $\mathcal{D}(\nu,\nu')=0$ if and only if $\nu=\nu'$. Moreover, it is also clear that we have
$$\mathcal{D}(\nu,\nu'')\le \mathcal{D}(\nu,\nu')+\mathcal{D}(\nu',\nu'')\cdot$$
So that $\mathcal{D}$ could be seen as a metric on $\fonc$  if $\mathcal{D}(\nu,\nu')<\pinf$ for any $\nu$ and $\nu'$, which is not the case. Nevertheless if $\mathcal{D}(\nu,\nu_{k})$ converges to $0$, when $k\to\pinf$, this tells us that $(\nu_{k})$ is ``converging'' to $\nu$.

Let $\ds Y=y_{1}\cdots y_{n}$,  $Y_{0}=\emptyset$ and, for any $1\le k \le p$, $Y_{k}=y_{1}\cdots y_{k}$. Since 
$$\nu(Y)=\prod_{k=1}^{p}\nu_{Y_{i-1}}(Y_{i}),$$
we conclude that for any $X\in\mathcal{C}$ and any $Y$we have
$$\left|\log\frac{\nu_{X}(Y)}{\nu'_{X}(Y)}\right|\le |Y|\mathcal{D}(\nu,\nu'),$$

This also implies, for any $XY\in \mathcal{C}_{n+p}$, with $|Y|=n$, that we have
$$\left|\nu_{X}(Y)-\nu'_{X}(Y)\right|
\le \nu_{X}(Y)\left(e^{n\mathcal{D}(\nu,\nu')}-1\right)\cdot$$
Let $g_{n}(X):=\scaln{\nu_{X}}{\log s_{X}}{n}$ and $g_{n}'(X):=\scaln{\nu_{X}'}{\log s_{X}'}{n}$. 
Assume that $\nu$ and $\nu'$ are both finite measures, and that $s_{X}\ge 1$, then we have
$$g_{n}(X)-g_{n}'(X)=\scaln{\left(1-\frac{\nu_{X}'}{\nu_{X}}\right)\nu_{X}}{\log s_{X}}{n}+\scaln{\nu_{X}'}{\log\frac{s_{X}'}{s_{X}}}{n},$$
which leads to
$$\left|g_{n}(X)-g_{n}'(X)\right|\le \left(e^{n\mathcal{D}(\nu,\nu')}-1\right)g_{n}(X)+n\mathcal{D}(s,s')\le \left(e^{n\mathcal{D}(\nu,\nu')}n\mathcal{D}(\nu,\nu')g_{n}(X)+n\mathcal{D}(s,s')\right),
$$
so that 
\begin{eqnarray}
\frac{1}{n}\left|g_{n}(X)-g_{n}'(X)\right| \le  e^{n\mathcal{D}(\nu,\nu')}\mathcal{D}(\nu,\nu')g_{n}(X)+\mathcal{D}(s,s')\label{conv1}
\cdot
\end{eqnarray}

\begin{Th}\label{thcont}
Let $\nu$ and $s$ be two elements of $\fonc$, and let $(\nu_{k})$ and $(s_{k})$ be two sequences of elements in $\fonc$.

For any $n\in \N$, set  $g_{n}=\scaln{\nu_{\emptyset}}{\log s_{\emptyset}}{n}$ and $g_{k,n}=\scaln{\nu_{k,{\emptyset}}}{\log s_{k,{\emptyset}}}{n}$.

Assume that 
\begin{itemize}
\item[--] $\nu$ and $\nu_{n}$ are finite measures.
\item[--] all those elements of $\fonc$ verify (ASI), with a same sequence $(\beta_{n})$,
\item[--] there exists a sequence $(\gamma_{n})$, compatible with $(\beta_{n})$, such that for any $s'\in \{s,\,s_{n}\}_{n\in \N}$ and any $XY\in\mathcal{C}$
$$\left|\log s'_{X}(Y)\right|\le \gamma_{|Y|},$$
\item[--] $\ds\lim_{k\to\pinf}\mathcal{D}(\nu,\nu_{k})=\lim_{k\to\pinf}\mathcal{D}(s,s_{k})=0$.
\end{itemize}
Then we have
$$\lim_{k\to\pinf}\limsup_{n}\left|\frac{g_{n,k}}{n}-\frac{g_{n}}{n}\right|=0,$$
And in particular we have
$$\lim_{k\to\pinf}\liminf_{n}\frac{g_{n,k}}{n}=\liminf_{n}\frac{1}{n}g_{n}\qquad\mbox{and}\qquad\lim_{k\to\pinf}\limsup_{n}\frac{g_{n,k}}{n}=\limsup_{n}\frac{1}{n}g_{n}\cdot$$
Moreover, the conclusion are still valid for any $\tilde{s}$ and $({\tilde{s_n}})$ such that $\tilde{s}\sim s$ and $\tilde{s}_n\sim s_n$ for any $n$.
\end{Th}
\proof By hypothesis we are in position to apply proposition \ref{sameasymptoticvalues} with $\nu$ and $s$, and $\nu_{n}$ and $s_{n}$ and we will use the notation therein.

Let $(A_{n})$ be a sequence of cylinders, with $A_{n}\in\mathcal{C}_{n}$ for any $n$, and let $p\le n$, with $q\in \N$ such that $qp<n\le (q+1)p$. Using (\ref{massexp4}) we get
$$g_{n}(\emptyset):=g_{n}=g_{pq+r}=g_{pq}+\sum_{|Y|=pq}\nu(Y)g_{r}(Y),$$
and also
$$g_{pq}=g_{p(q-1)}+\sum_{|Y|=p(q-1)}\nu(Y)g_{p}(Y),$$
so that
$$g_{pq}=\sum_{i=0}^{q-1}\sum_{|Y|=ip}\nu(Y)g_{p}(Y)\cdot$$
From this we deduce that 
$$g_{n}-\sum_{i=1}^{q-1}g_{p}(A_{ip})=\sum_{i=0}^{q-1}\sum_{|Y|=ip}\nu(Y)(g_{p}(Y)-g_{p}(A_{ik}))+\sum_{|Y|=pq}\nu(Y)g_{r}(Y),$$
and thus

$$
\left| \frac{g_{n}}{n}-\frac{1}{n}\sum_{i=1}^{q-1}g_{p}(A_{ip})\right|\le \frac{\alpha_{p}(\nu,s)}{p}+\frac{\gamma_{p}}{n}\le C\frac{c(\gamma_{p},\beta_{p})}{p} +\frac{\gamma_{p}}{n}\cdot
$$

The same is clearly also true for $\nu_{k}$ and $s_{k}$ so that we have 
$$
\begin{array}{rcl}
\ds\left| \frac{g_{n}}{n}-\frac{g_{k,n}}{n}\right| & \le & \ds\left| \frac{g_{n}}{n}-\frac{1}{n}\sum_{i=1}^{q-1}g_{p}(A_{ip})\right|+\left| \frac{g_{k,n}}{n}-\frac{1}{n}\sum_{i=1}^{q-1}g_{k,p}(A_{ip})\right|+\frac{1}{n}\left|\sum_{i=1}^{q-1}(g_{p}(A_{ip})-g_{k,p}(A_{ip}))\right|\\
& \le & \ds2\frac{\gamma_{p}}{n}+C\frac{c(\gamma_{p},\beta_{p})}{p}+\frac{1}{n}\left|\sum_{i=1}^{q-1}(g_{p}(A_{ip})-g_{k,p}(A_{ip}))\right|
\end{array}
$$
Using (\ref{conv1}) and the fact that one may assume that $\mathcal{D}(\nu,\nu_{k})\le 1$, we get
$$
\begin{array}{rcl}
\ds\left| \frac{g_{n}}{n}-\frac{g_{k,n}}{n}\right| & \le & \ds 2\frac{\gamma_{p}}{n}+C\frac{c(\gamma_{p},\beta_{p})}{p}+\frac{1}{n}\sum_{i=1}^{q-1}p\left(
e^{p\mathcal{D}(\nu,\nu_{k})}\mathcal{D}(\nu,\nu_{k})g_{p}(A_{ip})+\mathcal{D}(s,s_{k})\right)\\
 & \le & \ds 2\frac{\gamma_{p}}{n}+C\frac{c(\gamma_{p},\beta_{p})}{p}+e^{p}\gamma_{p}\mathcal{D}(\nu,\nu_{k})+\mathcal{D}(s,s_{k})
\end{array}
$$
Let $\varepsilon>0$.

Fix $p$ big enough to ensure that $\ds C\frac{c(\gamma_{p},\beta_{p})}{p}< \frac{\varepsilon}{3}$ and let $n_{0}\in\N$ such that $2\ds\frac{\gamma_{p}}{n_{0}}<\frac{\varepsilon}{3}$. 

Finally let $k_{0}$ such that for any $k\ge k_{0}$ we have $\ds e^{p}\gamma_{p}\mathcal{D}(\nu,\nu_{k})+\mathcal{D}(s,s_{k})<\frac{\varepsilon}{3}$.

Then for any $k\ge k_{0}$ and any $n\ge n_{0}$ we have : $\ds \left| \frac{g_{n}}{n}-\frac{g_{k,n}}{n}\right|<\varepsilon$, from where the result follows.
\qed

\subsection{Almost sure convergence}
In that section we work with a probability measure $\mu$. We are interested in  almost sure asymptotic values of sequences of the form : $Z_{n}(x):=\log s(X_{n}(x))$, where $s\in \fonc$, with $s\ge 1$, and for any $\mathtt{x}\in\K$, $X_{n}(\mathtt{x})$ is the unique element of $\mathcal{C}_{n}$ containing $\mathtt{x}$.
We will use the following strong law of large numbers version, see for instance \cite{HH} :

\begin{Th}\label{SLN}
Let $\left\{ Y_{n}, n \ge 1\right\}$ be a sequence of random variables in a probability space. Let $\left\{ \mathcal{F}_{n}, n \ge 1\right\}$ be an increasing sequence of $\sigma$-fields with $Y_{n}$ measurable with respect to $\mathcal{F}_{n}$, for each $n$. Assume that $(Y_{n})$ is uniformly bounded. Then almost-surely
$$
\lim_{n\to\pinf}\frac{1}{n} \sum_{i=1}^{n}\left(Y_{i}-\E_{\mu}\left[Y_{i}  |\mathcal{F}_{i-1}\right]\right)=0\cdot
$$
\end{Th}

\begin{Th}\label{thasi}
Let $\mu$ be a probability measure and let $s\in\fonc$.

Assume that $\mu$ and $s$ fulfilled (ASI) with sequence $(\beta_{n})$, that there exists a sequence $(\gamma_{n})$, compatible with $(\beta_{n})$ and such that  :
$$
\begin{array}{rc}
  \qquad & \forall XY\in\mathcal{C}\quad \left|\log s_{X}(Y)\right|\le \gamma_{|Y|}\\
\end{array}
$$
Then, for any $\tilde{s}\sim s$, there is a set $\mathcal{U}\subset \K$ of full $\mu$-measure such that for any $\mathtt{a}\in U$ we have :
$$\lim_{n\to\pinf}\left(\frac{1}{n}\sum_{|Y|=n}\mu(Y)\log \tilde{s}(Y)- \frac{1}{n}\log \tilde{s}(X_{n}(\mathtt{a}))\right)= 0\cdot$$
\end{Th}
\proof Given $n\in \N$, let $Z_{n}:=\frac{1}{n}\log s(X_{n})$ and, for any $p\in \N^*$,  $Y_{n,p}:=\ds\frac{1}{p}\left(Z_{np}-Z_{(n-1)p}\right)$. Let also $\mathcal{F}_{n}$ be the $\sigma$-algebra generated by cylinders in $\mathcal{C}_{n}$.
The idea, once again taking from \cite{Batakis5}, is to apply theorem \ref{SLN} on $(Y_{n,p})_{n}$, with respect to $(\mathcal{F}_{np})_{n}$, and to prove that the conditional expectations are asymptotically, on $n$ and $p$, independent of $\mathtt{x}\in\K$.

Fix $p>0$ and $k>0$ in $\N$. We have :
$$\E_{\mu}\left[Y_{k+1,p}  |\mathcal{F}_{kp}\right]=\frac{1}{p}\sum_{|Y|=p}\frac{\mu(X_{kp}Y)}{\mu(X_{kp})}\log \frac{s(X_{kp}Y)}{s(X_{kp})}=\frac{1}{p}g_{p}(X_{kp})\cdot$$
Recall that by (\ref{massexp4}), for any integer $k$, we have
$$g_{k+p}(X)=g_k(X)+\sum_{|Y|=k}\nu_{X}(Y)g_{p}(XY)\cdot$$
For any fixed cylinders $Z$ and $Y$ in $\mathcal{C}_{k}$ we have : $\ds| g_{p}(XY)-g_{p}(XZ)|\le \alpha_{p}$. From where we deduce that
$$\left|g_{k+p}(X)-g_{k}(X)-g_{p}(XZ)\right|=\left|\sum_{|Y|=k}\nu_{X}(Y)\left(g_{p}(XY)-g_{p}(XZ)\right)\right|\le \alpha_{p}.$$
Applying this with $XZ=X_{kp}:=X_{p}Y_{(k-1)p}$, we end up with :
$$\left|g_{kp}(X_{p})-g_{(k-1)p}(X_{p})-g_{p}(X_{kp})\right|\le \alpha_{p}.$$
Summing from $1$ to $q$ we get
$$\left|g_{qp}(X_{p})-g_{0}(X_{p})-\sum_{k=1}^qg_{p}(X_{kp})\right|\le q\alpha_{p}\cdot$$
This finally leads to
\begin{eqnarray}\label{esperance1}
\left|\frac{1}{q}\sum_{k=1}^q \E_{\mu}\left[Y_{k+1,p}  |\mathcal{F}_{kp}\right]-\frac{1}{qp}g_{qp}(X_{p})\right|\le \frac{\alpha_{p}}{p}\cdot
\end{eqnarray}
Let $n> p$ be two integers, and let $q\in \N$ be such that : $qp< n\le (q+1)p$. We have 
$$Z_{n}=Z_{(q+1)p}+Z_{n}-Z_{(q+1)p}=p\sum_{k=1}^q Y_{k+1,p}+\log \frac{s(X_{n})}{s(X_{(q+1)p})},$$
so that
$$
\frac{1}{n}Z_{n}=\frac{pq}{n}\frac{1}{q}\sum_{k=1}^q Y_{k+1,p}-\frac{1}{n}\log s_{X_{n}}(T^nX_{(q+1)p}),
$$
Using the hypothesis on $s$ we get
\begin{eqnarray}\label{esperance2}
\left|\frac{1}{n}Z_{n}-\frac{pq}{n}\frac{1}{q}\sum_{k=1}^q Y_{k+1,p}\right|\le C\frac{\gamma_{p}}{n}\cdot
\end{eqnarray}

Let 
$$
\begin{array}{rcl}
\mbox{\ding{192}} & = & \ds\frac{1}{n}Z_{n}-\frac{pq}{n}\frac{1}{q}\sum_{k=1}^q Y_{k+1,p}\\
 & &\\
\mbox{\ding{193}}_{q} & = &\ds \frac{1}{q}\sum_{k=1}^q Y_{k+1,p}-\E_{\mu}\left[Y_{k+1,p}  |\mathcal{F}_{kp}\right],\\
 & &\\
\mbox{\ding{194}} & = & \ds\E_{\mu}\left[Y_{k+1,p}  |\mathcal{F}_{kp}\right]-\frac{1}{qp}g_{qp}(X_{p}),\\
 & &\\
\mbox{\ding{195}} & = & \ds\frac{1}{n}\left(g_{qp}(X_{p})-g_{n}(X_{p})\right),\\
\end{array}
$$
We obviously have 
$$\frac{1}{n}Z_{n}-\frac{1}{n}g_{n}(X_{p})=\mbox{\ding{192}}+\frac{pq}{n}\mbox{\ding{193}}_{q}+\frac{pq}{n}\mbox{\ding{194}}+\mbox{\ding{195}}$$
So that, using (\ref{esperance1},\ref{esperance2}), we get
$$\left|\frac{1}{n}Z_{n}-\frac{1}{n}g_{n}(X_{p})\right|\le C\frac{\gamma_{p}}{n}+|\mbox{\ding{193}}_{q}|+\frac{\alpha_{p}}{p}+|\mbox{\ding{195}}|$$
Concerning \ding{195}, we use once again (\ref{massexp4}) with $r=n-pq\le p$
$$
g_{n}(X_{p})=g_{pq}(X_{p})+\sum_{|Y|=pq}\nu_{X_{p}}(Y)g_{r}(X_{p}Y)\cdot
$$
From where we deduce that
$$
|\mbox{\ding{195}}|=\frac{1}{n}\left|\sum_{|Y|=pq}\nu_{X_{p}}(Y)g_{r}(X_{p}Y)\right|\le C\frac{\gamma_{p}}{n},
$$
so that we have
$$
\left|\frac{1}{n}Z_{n}-\frac{1}{n}g_{n}(X_{p})\right|\le \frac{\alpha_{p}}{p}+C\frac{\gamma_{p}}{n}+|\mbox{\ding{193}}_{q}|\cdot
$$
Using (\ref{ineq2}) we get
\begin{eqnarray}\label{esperance3}
\begin{array}{rcl}
\ds\left|\frac{1}{n}Z_{n}-\frac{1}{n}g_{n}\right| & \le &\ds \frac{1}{n}\left|g_{n}-g_{n}(X_{p})\right|+\frac{\alpha_{p}}{p}+C\frac{\gamma_{p}}{n}+|\mbox{\ding{193}}_{q}|\\
 & \le & \ds C\frac{\alpha_{n}}{n}p+\frac{\alpha_{p}}{p}+C\frac{\gamma_{p}}{n}+|\mbox{\ding{193}}_{q}|\\
\end{array}
\end{eqnarray}
Note that for any $n$ we have by hypothesis that 
$$\ds \left|Y_{k,p}\right|= \left|\log\frac{s(X_{kp})}{s(X_{(k-1)p})}\right|= \left|\log s_{X_{(k-1)p}}(X_{kp})\right|\le \gamma_{p}\cdot$$
We may thus apply theorem \ref{SLN} with $(Y_{n,p})_{n}$, to conclude that $p$ being fixed, $|\mbox{\ding{193}}_{q}|$ $\mu$-almost surely converges towards $0$ as $q\to\pinf$, and we can easily conclude that there exists a set of full measure, denoted $\mathcal{U}$, such that it occurs for any $p$.

Let $\mathtt{x}\in \mathcal{U}$ and $\varepsilon>0$. Choose $p$ big enough to have $\ds\frac{\alpha_{p}}{p}<\frac{\varepsilon}{2}$. Remembering that $q$ goes to $\pinf$ with $n$, let $n_{0}\in \N$ such that for any $n\ge n_{0}$ we have 
$$\ds C\frac{\alpha_{n}}{n}p+C\frac{\gamma_{p}}{n}+|\mbox{\ding{193}}_{q}|<\frac{\varepsilon}{2}\cdot$$
Then by (\ref{esperance3}) we have
$$\left|\frac{1}{n}Z_{n}-\frac{1}{n}g_{n}\right| <\varepsilon,$$
and we can conlcude that the sequence $\ds\left(\frac{1}{n}Z_{n}(\mathtt{\cdot})-\frac{1}{n}g_{n}\right)$ $\mu$-almost surely converges towards $0$.
Which is precisely
$$\lim_{n\to\pinf}\left( \frac{1}{n}\log s(X_{n}(\mathtt{\cdot}))-\frac{1}{n}\sum_{|Y|=n}\mu(Y)\log s(Y)\right)= 0\cdot$$
\qed

\section{Pressures, Entropies, Lyapunov exponents and dimensions}

In this section, using results from the previous one, we show how natural generalizations of notions coming from Thermodynamical Formalism are possible.

\subsection{Lyapunov exponents and Entropies}
Given a sequence $(u_{n})$ of real numbers, we let $\lsup{u_{n}}$ be its limsup and $\linf{u_{n}}$ its liminf.\\

For any positive integer $n$ and any  probability measure, $\mu$, on $\K$  we define 
$$H_{\mu,n}:=\ds-\frac{1}{n}\sum_{|X|=n} \mu(X)\log\mu (X)\qquad \mbox{and}\qquad\Chi_{\mu,n}:=-\frac{1}{n}\sum_{|X|=n}\mu (X)\log \diame{X}{d},
$$
where the sum is taken over all cylinders $X\in \mathcal{C}_{n}$, with the usual convention $0\ln(0)=0$.

Let then $\lsup{h_{\mu}}:=\lsup{H_{\mu,n}}$, $\linf{h_{\mu}}:=\linf{H_{\mu,n}}$, $\lsup{\Chi_{\mu}}:=\lsup{\Chi_{\mu,n}}$ and $\linf{\Chi_{\mu}}:=\linf{\Chi_{\mu,n}}$.

Given any $\mathtt{a}\in \K$ and any positive number $n$, we define 
$$\chi_{n}(\mathtt{a}):=-\frac{1}{n}\log |X_{n}(\mathtt{a})|_{\mathsf{d}}\quad\mbox{and}\quad h_{\mu,n}(\mathtt{a}):=-\frac{1}{n}\log \mu(X_{n}(\mathtt{a}))\cdot$$
Let as before note : $\lsup{h_{\mu}}(\mathtt{a}):=\lsup{h_{\mu,n}(\mathtt{a})}$, $\linf{h_{\mu}}(\mathtt{a}):=\linf{h_{\mu,n}}(\mathtt{a})$, $\lsup{\Chi}(\mathtt{a}):=\lsup{\Chi_{n}}(\mathtt{a})$ and $\linf{\Chi}(\mathtt{a}):=\linf{\Chi_{n}}(\mathtt{a})$.
Remember that $\mathcal{F}_{n}$ is the $\sigma$-algebra generated by cylinders of length less or equal to $n$.  We have 
$$\ds H_{\mu,n}=\E_{\mu}\left[h_{\mu,n}\,|\, \mathcal{F}_{n}\right]\qquad\mbox{and}\qquad \Chi_{\mu,n} =\E_{\mu}\left[\chi_{n}\,|\, \mathcal{F}_{n}\right]\cdot$$ 
Note also that with $s:X\mapsto \ds\frac{1}{\diame{X}{d}}$ we have :
$$H_{\mu,n}=\frac{1}{n}\scaln{\nu}{\log \frac{1}{\nu}}{n}\qquad\mbox{and}\qquad \Chi_{\mu,n}=\frac{1}{n}\scaln{\nu}{\log s}{n}\cdot$$
Applying proposition \ref{jensen} to that setting, remembering that $D_n=d_0\cdots d_{n-1}$ is the degre of $T^n$, we get
\begin{Prop}
For any $n$ $$0\le H_{\mu,n}\le\frac{1}{n}\log D_{n}\quad  \mbox{and}\quad H_{\mu,n}-\Chi_{\mu,n}\le \frac{1}{n}\log\mathsf{m}_{n}(|.|_{\mathsf{d}})\cdot $$
Moreover, if  $X_{m}$ and $X_{M}$ in $\mathcal{C}_{n}$ are such that :  $\forall X\in\mathcal{C}_{n}$ we have $\diame{X_{m}}{d}\le\diame{X}{d}\le\diame{X_{M}}{d}$
then 
$$-\frac{1}{n}\log \diame{X_{M}}{d}\le \Chi_{\mu,n}\le -\frac{1}{n}\log \diame{X_{m}}{d}$$

\end{Prop}

\begin{Th}\label{entropy}
Let $\mu$ be a probability measure that fulfilled (ASI) with sequence $(\beta_{n})$ and assume that there exists a sequence $(\gamma_{n})$, compatible with $(\beta_{n})$, and such that  :

$$
\begin{array}{rc}
  \qquad & \forall XY\in\mathcal{C}\quad \left|\log \mu_{X}(Y)\right|\le \gamma_{|Y|}
\end{array}
$$

Then there exists a set $\mathcal{U}$ of full $\mu$-measure such that for any $\mathtt{a}\in\mathcal{U}$
$$\lim_{n\to\pinf} \left(h_{\mu,n}(\mathtt{a})- H_{\mu,n}\right)=0\cdot$$
and in particular we have $\mu$-almost surely
$$\lsup{h_{\mu}}(\mathtt{\cdot})=\lsup{h_{\mu}}\qquad \mbox{and}\qquad \linf{h_{\mu}}(\mathtt{\cdot})=\linf{h_{\mu}}\cdot$$
Assume moreover that $s:X\mapsto \ds\frac{1}{\diame{X}{d}}$ fulfilled (ASI) and that there exist $0<\delta<1$ such that for any $XY\in\mathcal{C}$ we have 
\begin{eqnarray}\label{diameter}
\delta^{\,\gamma_{|Y|}}\le \frac{\diame{XY}{d}}{\diame{X}{d}}\le 1
\end{eqnarray}
 then there exists a set of full $\mu$-measure, $\mathcal{V}$, such that for any  $\mathtt{a}\in\mathcal{V}$
$$\lim_{n\to\pinf}\left( \Chi_{n}(\mathtt{a})- \Chi_{\mu,n}\right)=0,$$
and in particular we have $\mu$-almost surely
$$\linf{\Chi}=\linf{\Chi}(\mathtt{\cdot})\le \lsup{\Chi}(\mathtt{\cdot}) =\lsup{\Chi}\cdot$$
\end{Th}
\proof 
The hypothesis allow us to apply theorem \ref{thasi}, with $\nu(X)=\mu(X)$ and $s(X)=\ds\frac{1}{\mu (X)}$, to conclude that
$$\lim_{n\to\pinf} \left(h_{\mu,n}(\mathtt{a})- H_{\mu,n}\right)=0\cdot$$

Let now $s:X\mapsto \ds\frac{1}{\diame{X}{d}}$.\\
From the assumption $\ds \delta^{\gamma_{|Y|}}\le \frac{\diame{XY}{d}}{\diame{X}{d}}\le 1$ we deduce that $1\le\ds s_{X}(Y)\le \delta^{-\gamma_{|Y|}}$, we easily conclude that for any $XY\in \mathcal{C}$, with $|Y|=n$ we have
$$\left|\log s_{X}(Y)\right|\le -\gamma_n\log d\cdot$$
We may once again apply theorem \ref{thasi} to conclude
\qed
\subsection{Pressure of a measure and dimensions}

In that paragraph we assume uniform bounds for the diameters, which is a very natural hypothesis in our context : there exists $\delta\in ]0,1[$ such that for any $n\in \N$, any cylinder $XY$ with $|Y|=n$ we have
\begin{eqnarray}\label{diameter2}
\delta^{n}\le \frac{\diame{XY}{d}}{\diame{X}{d}}\le (1-\delta)^n\cdot
\end{eqnarray}

In particular we thus have : $0<-\log (1-\delta)\le\linf{\Chi}\le \lsup{\Chi}\le -\log \delta$

\begin{Def}\label{presure1}
For any real number $t$ and any integer $n$ let $P_{\mu,n}(t):=H_{\mu,n}-t\Chi_{\mu,n}$.

We define $\linf{P_{\mu}}(t)=\linf{P_{\mu,n}}(t)$ and $\lsup{P_{\mu}}(t)=\lsup{P_{\mu,n}}(t)$.
\end{Def}

\begin{Prop}\label{dimension}
Assume that $\diame{.}{d}$ satisifes  condition (\ref{diameter2}) and that $\mu$ and $\diame{.}{d}$ both fulfill (ASI) with sequence $(\beta_n)$ compatible with $(\gamma_n=n)$, then the function $t\mapsto \linf{P_{\mu}}(t)$ and $t\mapsto \lsup{P_{\mu}}(t)$ are stricly decreasing and bi-Lipschitz

We have :
\begin{enumerate}
\item $0\le\linf{P_{\mu}}(0)=\linf{h_{\mu}}\le \lsup{P_{\mu}}(0)=\lsup{h_{\mu}}\cdot$
\item $\forall t\in \R\qquad \linf{h_{\mu}}-t\log (1-\delta)\le \linf{P_{\mu}}(t)\le \lsup{P_{\mu}}(t)\le \lsup{h_{\mu}}-t\log \delta\cdot$
\item Let $HD(\mu)$ be the Hausdorff dimension of $\mu$, and $PD(\mu)$ its packing-dimension, then
$$ \linf{P_{\mu}}(HD(\mu))=0=\lsup{P_{\mu}}(PD(\mu))\cdot$$
\item For any $\mathtt{a}$ in a set of full $\mu$-measure we have :
$$\liminf \frac{H_{\mu,n}}{\chi_{\mu,n}}=\liminf \frac{h_{\mu,n}(\mathtt{a})}{\chi_{n}(\mathtt{a})}=HD(\mu)\quad\mbox{and}\quad \limsup \frac{H_{\mu,n}}{\chi_{\mu,n}}=\limsup \frac{h_{\mu,n}(\mathtt{a})}{\chi_{n}(\mathtt{a})}=PD(\mu)\cdot$$
\end{enumerate}
\end{Prop}

\proof
From (\ref{diameter2}) we easily get for any $t< t+h$ and any $n\in \N$: 
$$-h\log (1-\delta)+P_{\mu,n}(t+h)\le P_{\mu,n}(t)\le P_{\mu,n}(t+h)-h\log \delta\cdot$$
From where we deduce, with $P$ standing either for $\linf{\P_\mu}$  or $\lsup{\P_\mu}$, that we have 
$$0<-h\log (1-\delta)\le P(t)-P(t+h)\le -h\log \delta\cdot$$
Which tells us that $\linf{\P_\mu}$ and $\lsup{\P_\mu}$ are stricly decreasing and bi-Lipschitz.

Points (1) and (2) are obvious from definitions and hypothesis.

We turn to results about dimensions.

Hypothesis on $\mu$ and $\diame{.}{d}$ allow us to apply theorem \ref{thasi} with, for any $X\in\mathcal{C}$, $\nu(X)=\mu(X)$ and $s(X)=\ds\frac{\diame{X}{d}^t}{\mu(X)}$, as well as theorem \ref{entropy}. In particular, $\mu$-almost surely we have
$$
\begin{array}{rrcl}
$\mbox{(E)}$ & \ds\lim_{k\to\pinf} h_{\mu,n_{k}}(\mathtt{\cdot})=h & \Longleftrightarrow & \ds\lim_{k\to\pinf} H_{\mu,n_{k}}=h,\\
$\mbox{(L)}$ &\ds\lim_{k\to\pinf} \Chi_{n_{k}}(\mathtt{\cdot})=\chi & \Longleftrightarrow &\ds \lim_{k\to\pinf} \Chi_{n_{k}}=\chi,\\
$\mbox{(P)}$ &\ds\lim_{k\to\pinf}\,\, \frac{1}{n_{k}}\log \frac{\diame{X_{n_{k}}(\cdot)}{d}^t}{\mu(X_{n_{k}}(\cdot))}=p & \Longleftrightarrow & \lim_{k\to\pinf} H_{\mu,n_{k}}-t\Chi_{\mu,n_{k}}=p
\end{array}
$$
Since $\linf{P_{\mu}}$ and $\lsup{P_{\mu}}$ are stricly decreasing, continuous with positive and negative values, there exists a unique $(\linf{t},\lsup{t})\in\R^2$, with $0\le \linf{t}\le\lsup{t}\in\R$ such that : $\linf{P_{\mu}}(\linf{t})=0=\lsup{P_{\mu}}(\lsup{t})$.
\\

\noindent We focus on the result with liminf.

Take a point $\mathtt{a}\in\K$ such that (E), (L) and (P) are true.

Let $(n_{k})$ be a sequence of integers such that 
$$\liminf \frac{h_{\mu,n}(\mathtt{a})}{\Chi_{n}(\mathtt{a})}=\lim_{k\to\pinf} \frac{h_{\mu,n_{k}}(\mathtt{a})}{\Chi_{n_{k}}(\mathtt{a})}=D(\mathtt{a})\cdot$$
Extracting again if necessary one may suppose that 
$$\ds\lim_{k\to\pinf} h_{\mu,n_{k}}(\mathtt{a})=h\,\mbox{ and } \lim_{k\to\pinf} \Chi_{\mu,n_{k}}(\mathtt{a})=\chi\cdot$$ By hypothesis on $\diame{.}{d}$, we have $0<\chi<\pinf$, thus $D(a)=\frac{h}{\chi}$ and we have 
$$\linf{P_{\mu}}(D(a))\le h-D(a)\chi=0\cdot$$
From where we deduce that $D(a)\ge \linf{t}$.

Assume now that $n_{k}$ is such that 
$$\lim_{k\to\pinf} H_{\mu,n_{k}}-t\Chi_{\mu,n_{k}}=\linf{P_{\mu}}(t)\cdot$$
Pick a point $\mathtt{a}$ for which we can defined $h$ and $\chi$ as above. We have $\ds\frac{h}{\chi}\ge D(\mathtt{a})$ and also
$$\linf{P_{\mu}}(t)=h-t\chi=\chi\left(\frac{h}{\chi}-t\right)\ge \chi(D(a)-t)\cdot$$
We deduce that $\linf{P_{\mu}}(D(a))\ge 0$ so that $D(a)\le\linf{t}\cdot$

This tells us that $\mu$-almost surely $D$ is constant equal to $\linf{t}$ which, by classical results on Hausdorff dimension of measure, see for instance \cite{falconer2013fractal}, allows us to conclude that $\mu$-almost surely $$HD(\mu)=D(\cdot)=\linf{t}\cdot$$
\qed
\subsection{Continuity of dimension}
\begin{Th}\label{contdim}
Let $\mu$ be a probability measure on $\K$ and $(\mu_{k})$ be a sequence of probability measures.

Let $\mathsf{d}$ be a metric on $\K$ and $\mathsf{d}_{k}$ a sequence of metrics. For any $X\in\mathcal{C}$ we let $\Tilde{s}(X)=\diame{X}{d}$ and $\Tilde{s}_{k}(X)=\diame{X}{d_{k}}$.

Assume that 
\begin{itemize}
\item[--]
The metrics $\diame{\cdot}{d}$ and $\diame{\cdot}{d_k}$, $k\in \N$, verify condition (\ref{diameter2}) with a common constant $\delta$.
\item[--]
$\mu$, $s$, and for any $k$, $s_{k}$ and $\mu_{k}$, fulfilled (ASI) with a common sequence $(\beta_{n})$.
\item[--]
there exists a sequence $(\gamma_{n})$, compatible with $(\beta_{n})$, such that for any $XY\in \mathcal{C}$ and any $k$,
$$
\log\frac{\mu(X)}{\mu(XY)}\le \gamma_{|Y|} \quad \mbox{and} \quad \log\ds\frac{\mu_{k}(X)}{\mu_{k}(XY)}\le \gamma_{|Y|}\cdot
$$
\item[--] 
we have $\ds\lim_{k\to\pinf}\mathcal{D}(\mu,\mu_{k})=\lim_{k\to\pinf}\mathcal{D}(s,s_{k})=0$.
\item[--]
$(t_{k})$ is a real sequence converging towards $t$.
\end{itemize}
Then 
$$\lim_{k\to\pinf}\limsup_{n}\left((H_{\mu,n}-t\chi_{\mu,n})-(H_{\mu_{k},n}-t_{k}\chi_{\mu_{k},n})\right)=0\cdot$$
In particular 
$$\lim_{k\to\pinf} \linf{P}_{\mu_{k}}(t_{k})=\linf{P}_{\mu}(t)\quad\mbox{and}\quad \lim_{k\to\pinf} \lsup{P}_{\mu_{k}}(t_{k})=\lsup{P}_{\mu}(t),$$
$$\lim_{k\to\pinf}\conf{k}{HD}{}(\mu_{k})=HD(\mu)\quad\mbox{and}\quad \lim_{k\to\pinf} \conf{k}{PD}{}(\mu_{k})=PD(\mu)\footnote{We would like to emphasize on the fact that the metrics changing with $k$, those dimensions are also changing, since Hausdorff and packing dimensions rely on the metric used.}.$$
\end{Th}
\proof  For any $X\in\mathcal{C}$ let 
$$\ds p(X)=\frac{\diame{X}{d}^{t}}{\mu(X)}=\frac{s(X)^t}{\mu(X)}\qquad\mbox{and for any $k$}\qquad p_{k}(X)=\frac{\diame{X}{d_{k}}^{t_{k}}}{\mu_{k}(X)}=\frac{s_{k}(X)^{t_{k}}}{\mu_{k}(X)}\cdot$$
Let $C=\ds(1+\sup_{k}\{|t|, |t_{k}|\})$. For any $XY\in\mathcal{C}$, we have 
$$\left|\log p_{X}(Y)\right|\le|t|\log s_{X}(Y)+\log\mu_{X}(Y)\le C\gamma_{|Y|}\cdot$$
In the same way we check that $\ds\left|\log p_{k_{X}}(Y)\right|\le C\gamma_{|Y|}\cdot$ Thus one may take $\gamma'_{n}=C\gamma_{n}$.

We also have for any $XYZ\in\mathcal{C}$ and any $X'\in\mathcal{C}_{|X|}$
$$
\begin{array}{rcl}
\ds \left|\frac{p_{XY}(Z)}{p_{X'Y}(Z)}-1\right| & = &\ds \left|\frac{\mu_{X'Y}(Z)}{\mu_{XY}(Z)}\left(\frac{s_{XY}(Z)}{s_{X'Y}(Z)}\right)^t-1\right|\\
 & = & \ds \left|\left(\left(\frac{s_{XY}(Z)}{s_{X'Y}(Z)}\right)^t-1\right)\frac{\mu_{X'Y}(Z)}{\mu_{XY}(Z)}+\frac{\mu_{X'Y}(Z)}{\mu_{XY}(Z)}-1\right|\\
 & \le & \ds (1+\beta_{0})\left|\left(\left(\frac{s_{XY}(Z)}{s_{X'Y}(Z)}\right)^t-1\right)\right|+\left|\frac{\mu_{X'Y}(Z)}{\mu_{XY}(Z)}-1\right|\\
 & \le & \ds 2(1+\beta_{0})^{1+|t|}\beta_{|Y|}\\
 & \le & \ds C\beta_{|Y|}
\end{array}
$$
We may thus take $\beta'_{n}=C\beta_{n}$. It is easy to prove that $(\gamma'_{n})$ and $(\beta'_{n})$ are compatible. We also have that $\mathcal{D}(p,p_{k})\to 0$, so that we can  apply theorem \ref{thcont}. Using notations therein  we have :
$$\frac{1}{n}g_{n}=H_{\mu,n}-t\chi_{\mu,n}=P_{\mu,n}(t)\qquad\mbox{and}\qquad \frac{1}{n}g_{k,n}=H_{\mu_{k},n}-t_{k}\chi_{\mu_{k},n}=P_{\mu_{k},n}(t_{k}),$$
and theorem \ref{thcont} asserts that
$$\lim_{k\to\pinf}\limsup_{n}\left(P_{\mu,n}(t)-P_{\mu_{k},n}(t_{k})\right)=0\cdot$$
From this we conclude that 
$$\lim_{k\to\pinf} \linf{P}_{\mu_{k}}(t_{k})=\linf{P}_{\mu}(t)\quad\mbox{and}\quad \lim_{k\to\pinf} \lsup{P}_{\mu_{k}}(t_{k})=\lsup{P}_{\mu}(t)\cdot$$
With $t_{k}=\conf{k}{HD}{}(\mu_{k})$, and $\delta$ any limit value of $(t_{k})$ we conclude, using proposition \ref{dimension}, that 
$$0=\lim_{k\to\pinf} \linf{P}_{\mu_{k}}(\conf{k}{HD}{}(\mu_{k}))=\linf{P}_{\mu}(\delta),
$$
which, by proposition \ref{dimension}, tells us that $\delta=HD(\mu)$, so that $(\conf{k}{HD}{}(\mu_{k}))$ converges towards $HD(\mu)$.
\qed
\section{Harmonic measure of non-autonomous conformal repellers}
We apply in this last section, results of the previous ones to study harmonic measure of the complementary of non-autonomous conformal repellers.

In the first paragraph we study geometric aspects of a given non-autonomous conformal iterated function system $\Psi$, proving that diameter and harmonic measure fulfilled the desired conditions.

Then we deal with sequence of such iterated function systems proving continuity of Hausdorff and Packing dimensions of harmonic measures.

Let $\Psi=({\Psi_n})=\{\conf{n}{\psi}{i}:U=\D\to U=\D\;;\; i=1,...,d_n\}$ be a non-autonomous conformal iterated function system satisfying our assumptions with a constant $\eta$ controlling contraction and modulus of the annulus, and a fixed neighborhood  of $\D$, $V$. We note $\X$ its limit set.
In this last section, we adopt the following notations:
\begin{itemize}
\item For any symbolic cylinder $X=a_{1}\cdots a_{n}$ we note $\psi_{X}$ the contraction $\conf{1}{\psi}{a_{1}}\circ \cdots\circ\conf{n}{\psi}{a_{n}}$, and $\phi_{X}=\conf{n}{\varphi}{a_{n}}\circ \cdots\circ\conf{1}{\varphi}{a_{1}}$ its inverse.
\item If $X$ is a symbolic cylinder, we also denote its geometric realization by $X$. That is, $X=\psi_{X}(\D)\cap \X$.
\item If $E\subset V$ and $X$ is a cylinder, then $E_{X}:=\psi_{X}(E)$. In particular, we denote $\D_{X}=\Psi_{X}(\D)$ and $\D_{\gamma X}=\Psi_{X}(\D_{\gamma})$, where $\D_{\gamma}=\{|z|\le \gamma\}$.
\item Given a Green domain $\Omega\subset \C$, a point $x\in\Omega$ and a Borel set $F$, we will denote $\omega(x,F,\Omega)$ the harmonic measure assigned to the domain $\Omega$ of $F$ at $x$.
\item In particular we will note $\omega:=\omega(\infty, \cdot, \hat{\C}\setminus \X)$, and for any $k$, $\omega_{k}:=\omega(\infty, \cdot, \hat{\C}\setminus \X_{k})$, where $\X_{k}=T^k(\X)$.
\end{itemize}
\subsection{Controlling diameters of cylinders}

As one may expect we will use K\oe be distortion theorem, and more precisely its quantified versions, to cope with diameters.

\begin{Prop}\label{koebe}
For any $0<\delta <1$ fixed, there exists $C>1$, such that for any conformal map, $f:\D\to \C$, any compact set $\K\subset\{|z|\le 1-\delta\}$ and any $(x,y)\in \K\times \K$
\begin{eqnarray}\label{kbe1}
 (1-C\diam (\K))\le \frac{|f'(x)|}{|f'(y)| }\le (1+C\diam (\K))
\end{eqnarray}

and 

\begin{eqnarray}\label{kbe2}
(1-C\diam (\K))|f'(x)|\le \frac{\diam f(\K)}{\diam \K }\le |f'(x)|(1+C\diam (\K))
\end{eqnarray}
\end{Prop}
\proof
First note that relation (\ref{kbe1}) and (\ref{kbe2}) are relevant only if $\diam (\K)$ is small. Otherwise it is nothing else than the standard Koebe distortion theorem. And we assume from now on that $\diam(\K)\le \frac{\delta}{2}$.

Let $x$ and $y$ be two points in $\K$. We apply corollary 1.5 in \cite{Pom1}, with $\ds z=\frac{x-y}{1-\bar{y}x}$, we have
$$
\frac{|f'(x)|}{|f'(y)|}\le \left(\frac{1+|z|}{1-|z|}\right)^3,
$$
Since $\K\subset \{|z|\le 1-\delta\}$ we have $|z|\le \frac{1}{\delta}|x-y|\le \frac{1}{\delta}\diam \K$. We also have $|z|<1$ and we can deduce that we have :
\begin{eqnarray}\label{KiF1}
\frac{|f'(x)|}{|f'(y)|}\le 1+C\diam (\K)\cdot
\end{eqnarray}
By exchanging $x$ and $y$ we also get the opposite inequality.

Let $Z=1-x\bar{y}$. Theorem 1 in \cite{KiFan} asserts that we have
\begin{eqnarray}\label{KiF2}
\frac{|Z|}{1-|x|^2}\frac{|Z|-|x-y|}{|Z|+|x-y|}
\le \frac{|f(x)-f(y)|}{|f'(x)||x-y|}\le \frac{|Z|}{1-|x|^2}\frac{|Z|+|x-y|}{|Z|-|x-y|}
\end{eqnarray}
Since $Z=1-|x|^2 -x(\widebar{y-x})$ we have 
$$1-|x|^2-|x-y|\le|Z|\le1-|x|^2+|x-y|$$
Which leads to
$$
1-\frac{|x-y|}{1-|x|^2}\le\frac{|Z|}{1-|x|^2}\le 1+\frac{|x-y|}{1-|x|^2}
.$$
But $|x|\le 1-\delta$, so $\delta\le 1-|x|^2\le 1$, and $|x-y|\le \diam \K$ which gives
\begin{eqnarray}\label{KiF3}
1-\frac{\diam\K}{\delta}\le\frac{|Z|}{1-|x|^2}\le 1+\frac{\diam\K}{\delta}
\end{eqnarray}
We have assumed that $\diam \K\le \frac{\delta}{2}$ thus
$$\frac{|Z|+|x-y|}{|Z|-|x-y|}=1+2\frac{|x-y|}{|Z|-|x-y|}\le 1+2\frac{|x-y|}{1-|x|^2-2|x-y|}\le 1+ \frac{2}{\delta(\delta-1)}\diam \K,$$
and we prove the same way that
$$1-\frac{2}{\delta}\diam \K\le \frac{|Z|-|x-y|}{|Z|+|x-y|}\cdot$$
Using those two last inequalities and (\ref{KiF3}) in (\ref{KiF2}) we end up with
$$
\left(1-\frac{2}{\delta}\diam \K\right)^2\le \frac{|f(x)-f(y)|}{|f'(x)||x-y|}\le\left(1+ \frac{2}{\delta(\delta-1)}\diam \K\right)^2\cdot
$$
Using (\ref{KiF1}), with an adjustement of the constant,  we may change $|f'(x)|$ with $|f'(z)|$, for any $z\in\K$, so that for any $x,y$ and $z$ in $\K$ we have
$$
1-C\diam \K\le \frac{|f(x)-f(y)|}{|f'(z)||x-y|}\le 1+C\diam \K\cdot
$$
We easily conclude that for any $z\in \K$
$$
\left(1-C\diam \K\right)|f'(z)|\le \frac{\diam f(\K)}{\diam \K}\le \left(1+C\diam \K\right)|f'(z)|\cdot
$$

\qed

We use this result to prove that asserts that relevant functionals verify (ASI).

For any $X\in\mathcal{C}$ let $s(X):=\diam (X)$, where $\diam$ is the euclidean diameter of the geometric cylinder $X$, and let $\Tilde{s}:=|(\psi_X)'(0)|$.

\begin{Cor}\label{derivKoebe}
There exists $C>1$ and $q\in ]0,1[$such that for any  cylinder $XYZ\in \mathcal{C}$ we have
\begin{eqnarray}\label{kbeCyl}
1-Cq^{|Y|}&\le& \Tilde{s}_{XY}(Z)\left|\frac{(\psi_ {Y})'(0)}{(\psi_ {YZ})'(0)}\right|\le 1+Cq^{|Y|}\nonumber \\
&\quad\mbox{and}\quad&  \\
1-Cq^{|Y|}&\le& s_{XY}(Z)\left|\frac{\diam (Y)}{\diam(YZ)}\right|\le 1+Cq^{|Y|}\nonumber
\end{eqnarray}

We have $s\sim \Tilde{s}$ and they both verify (ASI) with $\beta_n=Cq^n$.

\end{Cor}

\proof
If we apply classical Koebe distortion estimates with $f=\psi_X$,  $\K=\X_{k}$ and $x=0$ we get 
$$\frac{1}{C}\le \frac{\diam (X)}{|(\psi_X)'(0)|}\le C\cdot$$
Thus $s$ and $\Tilde{s}$ are equivalent.

We have 
$$
\left|\frac{(\psi_{XYZ})'(0)}{(\psi_{XY})'(0)}\right|=
\left|\frac{(\psi_{X})'(\psi_{YZ}(0))}{(\psi_{X})'(\psi_{Y}(0))}\right|\left|\frac{(\psi_{YZ})'(0)}{(\psi_{Y})'(0)}\right|\cdot
$$
Since $\psi_{YZ}(0)$ and $\psi_{Y}(0)$ both lie in $Y$, by (\ref{kbe1}) we have
$$(1-C\diam (Y))\le \left|\frac{(\psi_{X})'(\psi_{YZ}(0))}{(\psi_{X})'(\psi_{Y}(0))}\right|\le (1+C\diam (Y)),$$
from where we deduce, using also the fact that by hypothesis we have $\diam (Y)\le Cq^{|Y|}$, that 
$$(1-Cq^{|Y|})\left|\frac{(\psi_{YZ})'(0)}{(\psi_{Y})'(0)}\right| \le\left|\frac{(\psi_{XYZ})'(0)}{(\psi_{XY})'(0)}\right|\le (1+Cq^{|Y|})\left|\frac{(\psi_{YZ})'(0)}{(\psi_{Y})'(0)}\right|\cdot$$
This can be rewritten
$$1-Cq^{|Y|}\le \Tilde{s}_{XY}(Z)\left|\frac{(\psi_{Y})'(0)}{(\psi_{YZ})'(0)}\right|\le 1+Cq^{|Y|}\cdot$$
The same is also true with $X'$ instead of $X$ and this leads to
$$\left|\frac{\Tilde{s}_{XY}(Z)}{\Tilde{s}_{X'Y}(Z)}-1\right|\le Cq^{|Y|}\cdot$$
One can apply relation (\ref{kbe2}) with $f=\psi_{X}$,  $\K=YZ$, and also $\K=Y$,  and $x\in YZ$, to get
$$\frac{1-C\diam (YZ)}{1+C\diam (Y)}\le\frac{\diam (XYZ)}{\diam (XY) }\frac{\diam (Y)}{\diam (YZ) }\le \frac{1+C\diam (YZ)}{1-C\diam (Y)}.$$
We have $\diam (YZ)<\diam (Y)$, and by hypothesis we know that $\diam (Y)\le C q^{|Y|}$. Adjusting the constant $C$,  we easily conclude that we have
$$1-C q^{|Y|}\le \frac{\diam (XYZ)}{\diam (XY) }\frac{\diam (Y)}{\diam (YZ) }\le 1+C q^{|Y|}.$$
Which is exactly 
$$1-C q^{|Y|}\le s_{XY}(Z)\frac{\diam (Y)}{\diam (YZ) }\le 1+C q^{|Y|}.$$
As a matter of fact this implies that for any other geometric cylinder $X'\subset \X$, with $|X'|=|X|$, that we have
$$\left|\frac{s_{XY}(Z)}{s_{X'Y}(Z)}-1\right|\le C q^{|Y|}.$$

\qed

\subsection{Controlling diameters under perturbation}
We say that two NACIFS $\Psi$ and $\Tilde{\Psi}$ are compatible if for any $n$, 
\begin{itemize}
\item
 the degree of $\Phi_{n}$ and the degree of $\Tilde{\Phi}_{n}$ are the same,
\item
for any $i$ and $j$, $\conf{n}{\psi}{i}(V)\cap \conf{n}{\Tilde{\psi}}{j}(V)\neq \emptyset$ if and only if $i=j$.
\end{itemize}
For two compatible NACIFS $\Psi$ and $\Tilde{\Psi}$ we define
 $$d(\Psi,\Tilde{\Psi})=\sup_{n}\max_{1\le i\le d_{n}}\{\|\conf{n}{\psi}{i}-\conf{n}{\Tilde{\psi}}{i}\|_{\infty}\},$$
 with $\|f\|_{\infty}=\ds \sup_{z\in \bar{\D}_{\gamma}}|f(z)|$.
\begin{Prop}\label{contdiameters}
Assume that $\Psi$ and $\Tilde{\Psi}$ are two compatible NACIFS and that $0<\eta<1$ is the constant controling their annulus and their contractions. Let $\alpha=\ds \frac{\log (1-\eta)}{\log(\eta(1-\eta))}$. For any symbolic cylinder $X$ let $s(X)=|\psi_{X}'(0)|$ and $\Tilde{s}(X)=|\Tilde{\psi}_{X}'(0)|$. There exists $C>0$ such that 
$$\mathcal{D}(s,\Tilde{s})\le Cd(\Psi,\Tilde{\Psi})^{\alpha}\cdot$$
\end{Prop}

\proof Let $A_{n}=sup_{|X|=n}\|\psi_{X}-\Tilde{\psi}_{X}\|_{\infty}$. Note that we have $A_{1}=d(\Psi,\Tilde{\Psi})$. 

Let $Xa$ be a cylinder with $|X|=n$ and $|a|=1$. For any $z\in\D_{\gamma}$ we have
$$
\begin{array}{rcl}
\left|\psi_{Xa}(z)-\Tilde{\psi}_{Xa}(z)\right| &\le & \left|\psi_{X}(\psi_{a}(z))-\psi_{X}(\Tilde{\psi}_{a}(z))\right| +  \left|\psi_{X}(\Tilde{\psi}_{a}(z))-\Tilde{\psi}_{X}(\Tilde{\psi}_{a}(z))\right| \\
 & \le & (1-\eta)^n|\psi_{a}(z)-\Tilde{\psi}_{a}(z)|+A_{n}\\
 & \le & (1-\eta)^nA_{1}+A_{n},
\end{array}
$$
from where we deduce that 
\begin{equation}\label{majAn}
A_{n}\le \frac{1}{\eta}A_{1}=\frac{1}{\eta}d(\Psi,\Tilde{\Psi})
\end{equation}
Using Cauchy representation theorem, for any cylinder $X$ with $|X|=k$  we have
$$
\left|\psi_{X}'(0)-\Tilde{\psi}_{X}'(0)\right|= \frac{1}{2\pi}\left|\int_{\partial \D}\frac{\psi_{X}(z)-\Tilde{\psi}_{X}(z)}{z^2}dz\right|\le A_{k}.
$$
Since by assumption $\eta^k\le\ds \left|\Tilde{\psi}_{X}'(0)\right|$, we deduce that
\begin{equation}\label{maj1}
\left|\frac{\psi_{X}'(0)}{\Tilde{\psi}_{X}'(0)}\right|\le 1+\frac{1}{\eta^{k+1}}d(\Psi,\Tilde{\Psi})\cdot
\end{equation}
Inequality which is true for any cylinder $X$ such that $|X|\le k$, and stays true if we exchange the r\^ole played by $\psi$ and $\Tilde{\psi}$.

Let now $X=YX_{k}$ be any cylinder of length $n$, with $X_{k}=\emptyset$ if $n\le k$, and $|X_{k}|=k$ otherwise.
 
By chain rule we have
$$\psi_{YX_{k}}'(0)=\psi_{Y}'(\psi_{X_{k}}(0))\psi_{X_{k}}'(0)\cdot$$
So that for any cylinder $a$ of length 1 we have
$$
\frac{s_{X}(a)}{\Tilde{s}_{X}(a)}=\left|\frac{\psi_{Y}'(\psi_{X_{k}a}(0))}{\psi_{Y}'(\psi_{X_{k}}(0))}
\right|\left|\frac{\Tilde{\psi}_{Y}'(\Tilde{\psi}_{X_{k}}(0))}{\Tilde{\psi}_{Y}'(\Tilde{\psi}_{X_{k}a}(0))}\right|
\left|\frac{\psi_{X_{k}a}'(0)}{\Tilde{\psi}_{X_{k}a}'(0)}\right|
\left|\frac{\Tilde{\psi}_{X_{k}}'(0)}{\psi_{X_{k}}'(0)}\right|
$$
From (\ref{kbe1}) in proposition \ref{koebe} there exists $C>0$ such that
$$
\left|\frac{\psi_{Y}'(\psi_{X_{k}a}(0))}{\psi_{Y}'(\psi_{X_{k}}(0))}
\right| \le 1+C\mbox{diam} X_{k}\quad \mbox{and}\quad \left|\frac{\Tilde{\psi}_{Y}'(\Tilde{\psi}_{X_{k}}(0))}{\Tilde{\psi}_{Y}'(\Tilde{\psi}_{X_{k}a}(0))}\right|\le 1+C\mbox{diam} \Tilde{X}_{k}\cdot
$$
And from (\ref{maj1})
$$
\left|\frac{\psi_{X_{k}a}'(0)}{\Tilde{\psi}_{X_{k}a}'(0)}\right|
\left|\frac{\Tilde{\psi}_{X_{k}}'(0)}{\psi_{X_{k}}'(0)}\right|\le \left(1+\frac{1}{\eta^{k+1}}d(\Psi,\Tilde{\Psi})\right)^2\cdot
$$
Using the fact that $\mbox{diam} X_{k}$ and $\mbox{diam} \Tilde{X}_{k}$ are both smaller than $(1-\eta)^k$, and adjusting the constant $C$ we end up with
\begin{equation}\label{maj2}
\frac{s_{X}(a)}{\Tilde{s}_{X}(a)}\le 1+C\left((1-\eta)^k+\frac{1}{\eta^k}d(\Psi,\Tilde{\Psi})\right)=1+\frac{C}{\eta^k}\left((\eta(1-\eta))^k+d(\Psi,\Tilde{\Psi})\right)\cdot
\end{equation}
Let $k$ be the integer part of $\ds\frac{\log (d(\Psi,\Tilde{\Psi}))}{\log (\eta(1-\eta))}+1$, then $ k\ge\ds\frac{\log (d(\Psi,\Tilde{\Psi}))}{\log (\eta(1-\eta))}$ and $(\eta(1-\eta))^k\le d(\Psi,\Tilde{\Psi})$, so that we have 
$$
\frac{s_{X}(a)}{\Tilde{s}_{X}(a)}\le 1+2C\frac{d(\Psi,\Tilde{\Psi})}{\eta^k}\cdot
$$
On the other side we have $k\le \ds\frac{\log (d(\Psi,\Tilde{\Psi}))}{\log (\eta(1-\eta))}+1$ so that
$$
\frac{1}{\eta^k}\le \frac{1}{\eta}\exp\left(\frac{\log (d(\Psi,\Tilde{\Psi}))}{\log (\eta(1-\eta))}\log\frac{1}{\eta}\right)=\frac{1}{\eta}\left(\frac{1}{d(\Psi,\Tilde{\Psi})}\right)^{\frac{\log \eta}{\log(\eta(1-\eta))}}\cdot
$$
Which leads, with $\alpha=\frac{\log(1-\eta)}{\log(\eta(1-\eta))}$, to
$$
\frac{s_{X}(a)}{\Tilde{s}_{X}(a)}\le 1+Cd(\Psi,\Tilde{\Psi})^{\alpha},
$$
and gives
$$
\left|\log \left(\frac{s_{X}(a)}{\Tilde{s}_{X}(a)}\right)\right|\le Cd(\Psi,\Tilde{\Psi})^{\alpha}\cdot
$$
This being true for any cylinder $X$ we conclude that
$$
\mathcal{D}(s,\Tilde{s})\le Cd(\Psi,\Tilde{\Psi})^{\alpha}\cdot
$$
\qed

\subsection{Controlling harmonic measure}
In this section we extend well known properties of harmonic measure of usual IFS to NACIFS. We have tried to give all the ideas without drowning the reader with the total amount of technicalities.

Let $\Omega$ be a domain in $\mathbb{R}^2$  and 
$F \subset \partial\Omega$. For $x\in \Omega$ we denote  
$\omega(x,F,\Omega)$ the harmonic measure of $F$ in $\Omega$ evaluated at $x$. 
If $\Omega^c$ is bounded, we denote  $\omega(F,\Omega)$ the harmonic measure of $F$ 
in $\Omega$ evaluated at infinity.

We will frequently make use of the following (see for instance \cite{Brelot}): If $\Omega\subset\Omega'$ are two domains, $F\subset \partial\Omega\cap\partial \Omega'$ and $x\in\Omega$, then 
$$\omega(x, F,\Omega)= \omega(x, F,\Omega')-\int_{z\in\partial\Omega\cap\Omega'} \omega(z, F,\Omega')\omega(x, dz,\Omega)$$

For $\alpha,\beta$ positive real functions  we will write $\alpha \stackrel{c}{\sim} \beta$  if
$\frac1c\alpha \leq \beta \leq c \alpha$ with $c>0$ constant.
Finally, with no loss of generality we will take $U=\D$ the unit disc.

Remember that for any symbolic cylinder $X=a_1\cdots a_n$ we note $\Psi_X=\conf{1}{\psi}{a_{1}}\circ\cdots \circ\conf{n}{\psi}{a_{n}}$, and that we also denote $X$ the geometric cylinder $\X\cap \Psi_X(\D)$.

In the following we will deal with two NACIFS $\Psi$ and $\tilde\Psi$ and we will denote $\tilde X=\tilde\X\cap\tilde\Psi_X(\D)$ the geometric cylinder for $\tilde\Psi$. 

In order to apply the previous machinery we need to prove the following result

\begin{Th}\label{compare}
Let $\Psi$ and $\tilde\Psi$ be two NACIFS, $\X$, $\tilde\X$ the corresponding limit (Cantor)  sets and $\omega$, $\tilde\omega$ the harmonic measures on these sets.

We have :
\begin{equation}\label{compareform}
\lim_{d(\Psi, \tilde\Psi)\to 0}\mathcal{D}(\omega,\tilde\omega)=0\cdot
\end{equation}
\end{Th}

The proof follows the ideas of \cite{Batakis5} and heavily relies on a  repeated use of the following well known results.
\begin{Prop}\label{key} (\cite{MV},\cite{Ca})
Let $\Omega$ be  a domain containing $\infty$ and let 
$A_1 \subset B_1 \subset A_2\subset B_2 \subset ... \subset A_n\subset B_n\subset \Omega$
be conformal discs such that the annuli $B_i\setminus A_i$ are contained in $\Omega$,
for $1\le i\le n$. If the modules of the annuli are uniformly bounded away from zero and if
$\infty\in  \Omega \setminus B_n$
then, for all pairs of positive harmonic functions 
$u$, $v$ vanishing on  $\partial \Omega \setminus A_1$ and for
all $x \in \Omega \setminus B_n$ we have 
\begin{equation}
\left|\frac{u(x)}{v(x)}
: \frac{u(\infty)}{v(\infty)}-1\right|\leq {K  {q^n}}
\end{equation}
where $q<1$ and $K$ are two constants that depend only
on the lower  bound of the modules of the annuli.
\end{Prop}

\begin{Lemme}\label{minoring}
Let $\Psi$ be a NACIFS, $\mathbb{X}$ the corresponding limit (Cantor)  set and $\omega$ the harmonic measures on this sets. There exists $c>0$ depending only on $\gamma>1$, small enough so that $\gamma\D\subset V$, such that $\omega(x,\X,\D_\gamma\setminus\X)\ge c$ whenever $|x|=1$.
\end{Lemme}
\proof The proof implies standard arguments concerning capacity and harmonic measure on Cantor type sets: we hence leave some details to the reader (see also \cite{Batakis}).

Let $\mu$ be the uniform probability measure on $\X$, that is every cylinder of the $n^{\rm th}$ generation is charged by $1/(d_1\cdot \ldots \cdot d_n)$ mass.
For $x\ne y$ in $\D$ the Green function satisfies $\ln\left(\frac{C_1}{|x-y|}\right)\le G(x,y)\le \ln\left(\frac{C_2}{|x-y|}\right)$, where $C_1, C_2$ are constants depending on $\eta$.

Therefore, for $y\in\X$,  
$$\sum_{k=1}^{\infty}(d_1\cdot\ldots\cdot d_k)^{-1}\ln(C_1/(1-\eta)^k)\le G\mu(y)\le \sum_{k=1}^{\infty}(d_1\cdot\ldots\cdot d_k)^{-1}\ln(C_2/\eta^k),$$
which implies $G\mu\stackrel{c}{\sim} 1$ on $\X$. It follows on the maximum principle that $G\mu\stackrel{c}{\sim}\omega_{\D\setminus\X}(\X)$ in $\D$.
On the other hand, clearly $G\mu(x)\ge c>0$ for $x\in \D$  
and hence the result. 
\qed

The following lemma is a first application of this result. It will subsequently be used several times to obtain "localization" results for harmonic measure.

\begin{Lemme}\label{zooming}
Let $\Psi$ be a NACIFS and $\X$ its limit Cantor set. Assume that $\gamma>1$ is small enough so that $\D_{\gamma}:=\{|z|<\gamma\}\subset V$. There exists $C>0$ depending only on $\gamma$, so that for any cylinder $X\in \mathcal{C}$, and any $x\in\partial \D_{X}$ we have 
$$1\le\frac{\omega(x,X,\hat{\C}\setminus \X)}{\omega(x,X,\D_{\gamma X}\setminus X)}\le C\cdot$$
\end{Lemme}

\proof
Let $\Gamma := \{|x|=\gamma\}$, $u: x\mapsto\omega(x,X,\hat{\C}\setminus \X)$, $\tilde{u}:x\mapsto \omega(x,X,\D_{\gamma X}\setminus X)$ and $h=u-\tilde{u}$. The function $h$ is harmonic on $\D_{\gamma X}\setminus X$, with $h=0$ on $X$ and $h=u$ on $\partial \D_{\gamma X}$ so that
$$h(x)=\int_{\partial \D_{\gamma X}}u(z)\omega(x,dz,\D_{\gamma X}\setminus \X)\cdot$$
Let $x_{\gamma}\in\partial \D_{X}$ be such that for any $x\in \partial \D_{X}$, $u(x)\le u(x_{\gamma})$, and $x_{m}\in\partial \D_{\gamma X}$ be such that for any
$x\in \partial \D_{\gamma X}$ we have $u(x)\le u(x_{m})$.

The function $x\mapsto u(x)-u(x_{\gamma})$ is harmonic outside $\D_{X}\cup \X$, non-positive on $\partial \D_{X}\cup \X\setminus X$, it is thus negative on $\C\setminus \D_{X}$ and in particular we have $u(x_{m})<u(x_{\gamma})$.

We deduce the following
$$\tilde{u}(x_{\gamma})=u(x_{\gamma})-\int_{\partial \D_{\gamma X}}u(z)\omega(x_{\gamma},dz,\D_{\gamma X}\setminus \X)\ge u(x_{\gamma})\left(1-\omega(x_{\gamma},\partial \D_{\gamma X},\D_{X}\setminus \X)\right)=u(x_{\gamma})\omega(x_{\gamma},X,\D_{\gamma X}\setminus X)\cdot$$
By conformal invariance, with $n=|X|$, we have
$$\omega(x_{\gamma},X,\D_{\gamma X}\setminus \X)=\omega(\Phi_{X}(x_{\gamma}),\X_{n},\D_{\gamma}\setminus \X_{n}),$$
with by definition $|\Phi_{X}(x_{\gamma})|=1$. Lemma \ref{minoring} leads to $\ds\tilde{u}(x_{\gamma})\ge u(x_{\eta})c,$ since $\X_{n}$ is the limit Cantor set associated to an NACIFS with the same constant than $\X$.

We conclude that we have : $\ds 1\le \frac{u(x_{\gamma})}{\tilde{u}(x_{\gamma})}\le C\cdot$

By Harnack inequality  on $\partial\D_{X}$ we conclude that $\ds 1\le \frac{u(x)}{\tilde{u}(x)}\le C$ is true for any $x\in\partial\D_{X}$, with a different $C$, still uniquely depending on $\gamma$. 
\qed

We also need to establish the following lemma that restrains the number of generations necessary to establish an estimate on distribution of harmonic  measure.

\begin{Lemme}\label{localbehaviour}
Let $\Psi$ be a NACIFS, $\X$ its limit Cantor set and for any $\ell\ge 1$ let  $\X_{\ell}=T^{\ell}(\X)$.

There exist $q<1$ and $K>0$ depending only on $\eta$ such that for any  $XYZ\in \mathcal{C}$ with $|X|=k$, $|Y|=n$ we have 
\begin{equation}\label{local}
 \left|\log \left(\frac{\omega\left(XYZ\right)}{\omega\left(XY\right)}\frac{\omega_{k}\left(Y\right)}{\omega_{k}\left(YZ\right)}\right)\right|<Kq^{n}
\end{equation}
\end{Lemme}

The proof of this lemma follows the lines of the one in \cite{Batakis5}, therefore we will only give the highlights and main ideas of proof to facilitate the reading. 

\proof
Let $u:x\mapsto \frac{\omega\left(x,XYZ,\hat{\C}\setminus\X\right)}{\omega(XYZ)}$ and $v:x\mapsto\frac{ \omega\left(x,XY,\hat{\C}\setminus\X\right)}{\omega(XY)}$. Note that $u(\infty)=v(\infty)=1$.

Assume that $Y=y_{1}\cdots y_{n}$ and for any $\ell\in [1,n]$ let $X_{0}=X$, $X_{\ell}=Xy_{1}\cdots y_{\ell}$, $A_{\ell}=\D_{X_{n-\ell}}$ and $B_{\ell}=\D_{\gamma X_{n-\ell}}$. 

We have $A_{0}\subset B_{0}\subset \cdots\subset A_{n-1}\subset B_{n-1}\subset \hat{\C}\setminus \X$, and each annulus $B_{i}\setminus A_{i}$ is conformally equal to $\D_{\gamma}\setminus \D$.

The function $u$ and $v$ are both harmonic in $\hat{\C}\setminus\X$, and they vanished on $\X\setminus XY$.

We are thus in position to apply proposition \ref{key} to conclude that for any $x\notin \X\cup \D_{\gamma Xy_{1}}$, and in particular for any $x\in\partial\D_{X}\cup\partial\D_{\gamma X}$

\begin{equation}\label{step1}
\left|\log \left(\frac{\omega\left(x,XYZ,\hat{\C}\setminus\X\right)}{\omega\left(x,XY,\hat{\C}\setminus\X\right)}\frac{\omega(XY)}{\omega(XYZ)}\right)\right|=\left|\log \left(\frac{u(x)}{v(x)}\frac{v(\infty)}{u(\infty)}\right)\right|\le Cq^n\cdot
\end{equation}
Let now $\tilde{u}:x\mapsto \frac{\omega(x,XYZ,\D_{\gamma X}\setminus X)}{\omega(XYZ)}$ and $\tilde{v}:x\mapsto \frac{\omega(x,XY,\D_{\gamma X}\setminus X)}{\omega(XY)}$.
We have for any $x\in \D_{\gamma X}\setminus X$
$$\tilde{u}(x)=u(x)-\int_{\partial\D_{\gamma X}}u(z)\omega(x,dz,\D_{\gamma X}\setminus X),$$
and also
$$\tilde{v}(x)=v(x)-\int_{\partial\D_{\gamma X}}v(z)\omega(x,dz,\D_{\gamma X}\setminus X).$$
We thus have
$$\tilde{u}(x)-\tilde{v}(x)=u(x)-v(x)-\int_{\partial\D_{\gamma X}}(u(z)-v(z))\omega(x,dz,\D_{\gamma X}\setminus X).$$
Which may be written :
$$\frac{\tilde{u}(x)}{\tilde{v}(x)}-1=\frac{v(x)}{\tilde{v}(x)}\left(\left(\frac{u(x)}{v(x)}-1\right)-\int_{\partial\D_{\gamma X}}\frac{v(z)}{v(x)}\left(\frac{u(z)}{v(z)}-1\right)\omega(x,dz,\D_{\gamma X}\setminus X)\right).$$
Since (\ref{step1}) is true for any $x \in\partial\D_{X}\cup\partial\D_{\gamma X}$, and by lemma \ref{zooming}, we conclude that there exist  $C>0$ and $q<1$ such that for any $x\in \partial \D_{X}$
$$\left|\frac{\tilde{u}(x)}{\tilde{v}(x)}-1\right|\le Cq^n\left(1+\int_{\partial\D_{\gamma X}}\frac{v(z)}{v(x)}\omega(x,dz,\D_{\gamma X}\setminus X)\right)=Cq^n\left(1+\frac{v(x)-\tilde{v}(x)}{v(x)}\right)\le 2Cq^n\cdot$$
From this we conclude that for any $x\in \partial\D_{X}$ we have 
\begin{equation}\label{step2}
\left|\log\left( \frac{\omega\left(x,XYZ,\D_{\gamma X}\setminus X\right)}{\omega\left(x,XY,\D_{\gamma X}\setminus X\right)}\frac{\omega(XY)}{\omega(XYZ)}\right)\right|=\left|\log \left(\frac{\tilde{u}(x)}{\tilde{v}(x)}
\right)\right|\le Cq^n\cdot
\end{equation}
By conformal invariance of harmonic measure, we conclude that for any $x\in\partial \D$
$$
\left|\log \left(\frac{\omega\left(x,YZ,\D_{\gamma}\setminus \X_{k}\right)}{\omega\left(x,Y,\D_{\gamma}\setminus \X_{k}\right)}\frac{\omega(XY)}{\omega(XYZ)}\right)\right| \le Cq^n\cdot
$$
Applying  (\ref{step1}) with $\X_{k}$ and $\omega_{k}$, this is possible since the NACIFS giving birth to $\X_{k}$ is comparable to the one giving birth to $\X$, we get for any $x\in\partial \D$
\begin{equation}\label{step4}
\left|\log\left(\frac{\omega\left(x,YZ,\D_{\gamma}\setminus \X_{k}\right)}{\omega\left(x,Y,\D_{\gamma}\setminus \X_{k}\right)}\frac{\omega_{k}(Y)}{\omega_{k}(YZ)}\right)\right|\le Cq^n\cdot
\end{equation}
Combining (\ref{step1},\ref{step2},\ref{step4}) we end up with
\begin{equation}\label{step5}
\left|\log\left( \frac{\omega(XYZ)}{\omega(XY)}\frac{\omega_{k}(Y)}{\omega_{k}(YZ)}\right)\right|\le Cq^n\cdot
\end{equation}
\qed

From this result we easily deduce that harmonic measure verifies the two assumptions needed to apply main results of the previous section : (ASI) and uniform control on the measure of cylinders.
\begin{Cor}\label{measharmASI}
Let $\Psi$ be a NACIF, $\X$ its limit set and $\omega$ the harmonic measure on $\C\setminus \X$ evaluated at $\infty$. Under our geometric assumptions on the NACIF there exists $q\in ]0,1[$ and $C>0$ such that for any geometric cylinders $XYZ\in \mathcal{C}$ and $X'YZ\in \mathcal{C}$, with $|X|=|X'|$ and $|Y|=n$ 
\begin{align}
    \left|\log\left(\frac{\omega_{XY}(Z)}{\omega_{X'Y}(Z)}\right)\right| & \le  Cq^n\label{harASI1}\\
    \left| \log\left(\omega_{X}(Y)\right)\right| & \le  Cn \label{harASI2}
\end{align}
\end{Cor}

\proof
Relation (\ref{step5}) asserts that for any $XYZ\in \mathcal{C}$ with $|X|=k$ and $|Y|=n$: $$\ds\left|\log\left( \omega_{XY}(Z)\frac{\omega_{k}(Y)}{\omega_{k}(YZ)}\right)\right|\le Cq^n,
$$ 
from where we get, for any $X'\in \mathcal{C}$ with $|X'|=|X|$
$$
\begin{array}{rcl}
  \ds\left|\log\left(\frac{\omega_{XY}(Z)}{\omega_{X'Y}(Z)}\right)\right| & = &  \ds \left|\log\left( \omega_{XY}(Z)\frac{\omega_{k}(Y)}{\omega_{k}(YZ)}\right)-\log\left( \omega_{X'Y}(Z)\frac{\omega_{k}(Y)}{\omega_{k}(YZ)}\right)\right|   \\[0.5cm]
     &  \le & \ds\left|\log\left( \omega_{XY}(Z)\frac{\omega_{k}(Y)}{\omega_{k}(YZ)}\right)\right|+\left|\log\left( \omega_{X'Y}(Z)\frac{\omega_{k}(Y)}{\omega_{k}(YZ)}\right)\right|\\[0.5cm]
     & \le & Cq^n
\end{array}
$$

From (\ref{step5}) we deduce that there exists $C>1$ such that  $$\ds\frac{1}{C}\omega_k(Z)\le\frac{\omega(XZ)}{\omega(X)}\le 1\cdot$$

Assume that $Z=z_0\cdots z_{n-1}$ and for any $i\in \segN{0}{n-1}$ let $Z_i=z_0\cdots z_{i}$ then we have :
$$\omega_k(Z)=\omega_{k}(z_0)\Prod_{i=1}^{n-1}\frac{\omega_k(Z_i)}{\omega_k(Z_{i-1})},$$

since for any $i\in\segN{1}{n-1}$ relation (\ref{step5}) tells us that : $\ds\frac{1}{C} \omega_{k+i}(z_{i})\le \frac{\omega_k(Z_i)}{\omega_k(Z_{i-1})}\le 1$, we conclude that 
$$
\frac{1}{C^n}\Prod_{i=1}^{n-1}\omega_{k+i}(z_i)\le\frac{\omega(XZ)}{\omega(X)}\le 1\cdot
$$
From Lemma \ref{minoring} we can deduce that given a NACIF $\Psi$ there exists a constant $\alpha>0$, depending only the geometric constants that constrain the geometry of the NACIF, such that for any $\ell\ge 0$ and any $a\in \mathcal{C}_{\ell}$ with $|a|=1$, we have $\omega_{\ell} (a)\ge \alpha$, we conclude that
$$
\left(\frac{\alpha}{C}\right)^n\le \frac{\omega(XZ)}{\omega(X)}\le 1,
$$
and the result follows.

\qed

We turn to a comparison of harmonic measures of the limit sets $X$ and $\tilde X$. 

\begin{Lemme}\label{xi} Let $\gamma>1$ such that $\D_\gamma\subset V$ appearing in the definition of the NACIF $\Psi$.
There exist $K>0$ and $\alpha>0$ such that for any $x\in \D\setminus \X$ :
$$
\omega(x,\X,\D_\gamma\setminus \X)\ge 1-K\dist(x,\X)^\alpha\cdot
$$
\end{Lemme}

\proof
Given any $x\in \D\setminus \X$ let $n(x)=\displaystyle\max \{k\in\N\,|\,\exists X\in \mathcal{C}_n\,,\,x\in \D_X\}$.
Note that $n(x)$ is such that, there exists $X\in\mathcal{C}_n$ with $x\in \D_X$ and $x\notin \D_{Xa}$ for any cylinder $Xa \in \mathcal{C}_{n+1}$.

Using open set and bounded contraction hypothesis, we deduce from K\oe be distortion theorem that there exist three constants $K'>0$ and $0<\beta'<\beta $ such that
\begin{equation}\label{estim_n}
\frac{1}{K'}\dist(x,\X)^\beta\le e^{-n} \le K'\dist(x,\X)^{\beta'}
\end{equation}
It is important to note that $K'$, $\beta'$ and $\beta$ only depend on the constant $\eta$ used to control the geometry of the NACIF.

Let $x\in \D_X$, with $X=a_1\cdots a_n$ and $|X|=n(x):=n$. Let $X_0=\X$ and, for any $1\le k\le n$, $X_k=a_1\cdots a_k$ and set $u_k:y\mapsto \omega(y, X_k,\D_{\gamma X_k}\setminus X_k)$. The function $u_k$ is harmonic on $\D_{\gamma X_k}\setminus X_k$, and equal $1$ on $X_k$. Our goal is to prove that 
$$
u_0(x)=\omega(x,\X,\D_\gamma\setminus \X)\ge 1-K\dist(x,\X)^\alpha\cdot
$$
First we note, using  Lemma \ref{minoring}, conformal invariance and  maximum principle that there exists $c>0$ such that for any $x\in \D$
\begin{equation}\label{minor1}
    u_k(x)\ge c.
\end{equation}
For any $1\le k\le n$ we have
$$
u_{k-1}(x)-u_k(x)=\int_{\partial\D_{\gamma X_k}}u_{k-1}(z)\omega (x,dz,\D_{\gamma X_k}\setminus X_k )\cdot
$$
Because of our geometric assumptions on the NACIF and Harnack inequalities, there exists a constant $\lambda>0$ such that for any $Y\in\mathcal{C}$, any continuous  function $h>0$, defined on $\D_Y$ and harmonic on $\D_Y\setminus Y$,  and any $y$ and $y'$ in $\D_Y\setminus Y$ : $\displaystyle \lambda h(y)\le h(y')$.

From this we deduce that 
$$
u_{k-1}(x)-u_k(x)\ge \lambda u_{k-1}(x)\omega(x,\partial \D_{\gamma X_k},D_{\gamma X_k}\setminus X_k)=\lambda u_{k-1}(x)(1-u_k(x)).
$$
Let $v_k=1-u_k$, we drop $x$ to lighten notations. 
\\Using (\ref{minor1}) we get : $\lambda c v_k\le \lambda u_kv_k\le v_k-v_{k-1}$ and
$$
c\lambda\le \frac{v_k-v_{k-1}}{v_k}\le \int_{v_{k-1}}^{v_k}\frac{1}{t}dt=\log (v_k)-\log (v_{k-1})\cdot
$$
Summing for $1\le k\le n$
$$
c\lambda n\le \log v_n-\log v_0,
$$
so that
$$
1-u_0=v_0\le e^{-c\lambda n}.
$$
We use (\ref{estim_n}) to conclude that : $\displaystyle u_0(x)=\omega(x,\X,\D_\gamma\setminus \X)\ge 1-\left(\frac{1}{K'}\dist(x,\X)^\beta\right)^{c\lambda}$. This is the desired result with $K=\displaystyle\frac{1}{K'^{c\lambda}}$ and $\alpha=c\lambda\beta$.
\qed

This corollary happens to be crucial in order to prove continuity of dimensions.

\begin{Cor}\label{continuity1}
Let $\varepsilon>0$, there exists $n\in \N$ such that for any $XY\in \mathcal{C}$ with $|Y|=n$ and any $x\in\D_{XY}$ 
$$\omega(x,X,\C\setminus \X)\ge 1-\varepsilon\cdot $$
\end{Cor}

\proof
Let $\varepsilon >0$ and $\eta>0$ given by Lemma \ref{xi}. Let $n$ big enough so that for any $k\in \N$ and any $Y\in \mathcal{C}$ with $|Y|=n$ we have $\diam \D_{\gamma Y} < \eta$.
Let $X \in \mathcal{C}$ and $x\in \partial D_{XY}$, with $|Y|=n$. We have 
$$\omega (x,X,\C\setminus \X)\ge \omega (x,X,\D_{\gamma X}\setminus \X).$$
By conformal invariance, with $k=|X|$, we have 
$$\omega (x,X,\D_{\gamma X}\setminus \X)=\omega_k(\varphi_X(x),\X_k,\D_{\gamma}\setminus \X_k)\cdot$$
Since $x\in \partial \D_{XY}$ we have $\varphi_X(x)\in \partial\D_Y$, so that $\dist (\varphi_X(x),\X_k)\le \diam \D_{\gamma Y}<\eta$ and Lemma \ref{xi} gives : $\omega_k(\varphi_X(x),\X_k,\D_\gamma\setminus \X_k)>1-\varepsilon$.

We have proved that for any $x\in \partial D_{XY}$ we have 
$$
\omega (x,X,\C\setminus \X)\ge \omega (x,X,\D_{\gamma XY}\setminus \X)=\omega_k(\varphi_X(x),\X_k,\D_\gamma\setminus \X_k)\ge 1-\varepsilon\cdot
$$
This inequality extends to $\D_{XY}$ by maximum principle.

\qed

\noindent{\em Proof of theorem \ref{compare}. }
Let $\varepsilon >0$ and $N_1\in \N$ big enough so that  $Cq^{N_1}$ in Corollary \ref{measharmASI} is strictly less than $\varepsilon$.

From Corollary \ref{measharmASI} we know that : $cq^{N_1}\le \omega_k(Z)$ for any $k\in\N$, any $Z\in T^k(\mathcal{C})$ with $|Z|\le N_1$. Let $\varepsilon' =cq^{N_1}\varepsilon$. For any $k\in\N$, any $Z\in T^k(\mathcal{C})$ with $|Z|\le N_1$ we have  
\begin{equation}\label{epsilon}
    \frac{\varepsilon'}{\omega(X)}\le \varepsilon.
\end{equation}

Let $N_2$ be given by Corollary \ref{continuity1} with $\varepsilon'>0$ so that for any $XY\in\mathcal{C}$ with $|Y|=N_2$ and any $x\in \partial \D_{XY}$ we have :  $\omega(x,X,\C\setminus \X)\ge 1-\varepsilon'$.

Let now $d(\Psi,\tilde\Psi)$ be as small as necessary to ensure that for any $k\in \N$, any $Y\in T^k(\mathcal{C})$ with $|Y|\le N_1+N_2$ we have $\tilde Y\subset \D_{Y}$. Where $\tilde Y$ is the geometric cylinder associated to $\tilde\Psi$ with coding $Y$.

Let  $Xa\in \mathcal{C}$ with  $|a|=1$ and $|X|=k$. We have
$$
\log\left(\frac{\omega_{X}(a)}{\tilde\omega_{X}(a)}\right)=\log\left(\frac{\omega_{k}(a)}{\tilde\omega_{k}(a)}\right)+\log\left(\frac{\omega_{X}(a)}{\tilde\omega_{X}(a)}\frac{\tilde\omega_k(a)}{\omega_k(a)}\right)\cdot
$$ 
If $|X|=k>N_1$ we have by Corollary \ref{measharmASI}

\begin{equation}\label{conteq1}
   \log\left(\frac{\omega_{X}(a)}{\tilde\omega_{X}(a)}\right)\le \log\left(\frac{\omega_{k}(a)}{\tilde\omega_{k}(a)}\right)+\varepsilon
\end{equation}

We are thus led  to prove that $\ds \log\left(\frac{\omega_{X}(a)}{\tilde\omega_{X}(a)}\right)$ is small for any cylinder
$X$ with $|X|\le N_1$.

Actually we are going to prove that  $\ds\log\left(\frac{\omega_k(X)}{\tilde\omega_k(X)}\right)$ is small for any cylinder $X$ with $|X|\le N_1$ and any $k$.

In order to lighten notations, we drop the letter $k$. Remember that our estimates about harmonic measure are only depending on the geometric constants in our hypothesis : Bounded contraction and Annulus condition.

Let $|X|=k\le N_1$, $F=\ds\bigcup_{|Z|=N_1+N_2}\D_{\gamma Z}$ and $\Omega=\C\setminus F$. Our choices imply that for any $Z$ with $|Z|=N_1+N_2$ we have $\tilde Z \subset \D_Z$ so that $\tilde\X\subset F$.

The functions $u:x\mapsto \omega(x,X,\C\setminus \X)$ and $v:x\mapsto \omega(x,\tilde X,\C\setminus \tilde\X)$ are harmonic on $\Omega$. 

Let $|Z|=N_1+N_2$ with $Z\subset X$, then from Corollary \ref{continuity1} we deduce, for any $x\in \partial_{\gamma Z}$, that $u(x)\ge 1-\varepsilon$ as well as $v(x)\ge 1-\varepsilon'$, which implies that : $|u(x)-v(x)|\le \varepsilon'$ for any $x\in \partial_{\gamma Z}$.

Let now $Z=YW$, with $|Y|=|X|$, $|Z|=N_1+N_2$ and $Z\not\subset X$. By Corollary \ref{continuity1} we have $\omega(x,Y,\C\setminus \X)>1-\varepsilon'$ and also $\omega(x,\tilde Y,\C\setminus \tilde\X)>1-\varepsilon'$ for any $x\in \partial \D_{\gamma Z}$. 

Since 
$$1=\ds\sum_{|V|=|X|}\omega(x,V,\C\setminus\X)=\ds\sum_{|V|=|X|}\omega(x,\tilde V,\C\setminus\tilde\X),$$ we deduce that $u(x)\le 1-\omega(x,Y,\C\setminus \X)\le \varepsilon'$ and $v(x)\le 1-\omega(x,\tilde Y,\C\setminus \tilde \X)\le \varepsilon'$. We thus have, for any $x\in \partial \D_{\gamma Z}$,
$|u(x)-v(x)|\le 2\varepsilon'$.

The functions $u$ and $v$ are both harmonic on $\Omega$ and for any $x\in \partial \Omega$ we have $|u(x)-v(x)|\le 2\varepsilon'$. By maximum principle we have $|u(x)-v(x)|\le 2\varepsilon'$ for any $x\in \Omega$ and in particular for $x=\infty$, so that $|\omega(X)-\tilde\omega (X)|\le 2\varepsilon'$. 

We deduce from (\ref{epsilon}) that for any cylinder $X$ with $|X|\le N_1$ we have : 
$$\ds\left|\frac{\tilde\omega(X)}{\omega(X)}-1\right|\le 2\varepsilon'\frac{1}{\omega(X)}\le 2\varepsilon\cdot$$

From where it is easy to deduce that we have for such cylinders :
$$
\log\left(\frac{\omega_X(a)}{\tilde\omega_X(a)}\right)=\log\left(\frac{\omega(Xa)}{\tilde\omega(Xa)}\right)+\log\left(\frac{\tilde\omega(X)}{\omega(X)}\right)\le 4\varepsilon
$$

With inequality (\ref{conteq1}) we can conclude that for any cylinder $X$ of any length we have
$$
\log\left(\frac{\omega_X(a)}{\tilde\omega_X(a)}\right)\le 5\varepsilon.
$$
Which tells us that $\mathcal{D}(\omega,\tilde\omega)\le 5\varepsilon$ and concludes the proof.

\qed
\subsection{Back to the complex plane}
We are now in position to give a proof of the main results of this paper about harmonic measure on the plane: Theorems \ref{Continuity} and \ref{thdim}

Let $\X$ be a non-autonomous Cantor set  associated with the conformal iterated systems $\Psi$ and let $\omega$ be the harmonic measures of its complementary in $\C$.

We defined on the symbolic space $\K$ the metric $d$. For any $\mathtt{a}=(a_{n})\in\K$ and  $\mathtt{b}=(b_{n})\in\K$ let  $d(\mathtt{a},\mathtt{b})=|x(\mathtt{a})-x(\mathtt{b})|$, which is the euclidian distance in $\C$ of the points encoded by $\mathtt{a}$ and $\mathtt{b}$. In that way, symbolic cylinders and geometric cylinders have the same diameters. 

From the discussions in paragraph 4.1 we easily get :
$$\frac{\diame{XY}{d}}{\diame{X}{d}}\sim \frac{|\Psi_{XY}'(0)|}{|\Psi_X'(0)|}=\frac{|\Psi_{X}'(\Psi_Y(0))|}{|\Psi_X'(0)|}|\Psi_Y'(0)|\sim |\Psi_Y'(0)|\cdot
$$
Hypothesis Bounded Contraction implies that there exists a constant $\eta\in ]0,1[$ such that : $(1-\eta)^{|Y|}\le |\Psi_Y'(0)|\le \eta^{|Y|}$.

And we can conclude that $d$ satisfies condition (\ref{diameter2}).

By Corollary \ref{derivKoebe} we also know that $s(X):=\diame{X}{d}$ satisfies $(ASI)$ with sequence $(\beta_n=Cq^n)$ with $0<q<1$.

Corollary \ref{measharmASI} expresses the fact that $\omega$ also fulfills (ASI) with sequence $(\beta_n=Cq^n)$, with $0<q<1$, with an adjustment of the constant $q$, as well as condition (\ref{massexpbis1}) with $\gamma_n\sim Cn$.

We are thus in position to apply Theorem \ref{dimension} with $\mu=\omega$ to conclude that Theorem \ref{thdim} is true.

Let now $(\X^k)_{k\in\N}$ be a sequence of non-autonomous Cantor sets  associated with the conformal iterated systems  $(\Psi^k)_{k\in \N}$, and Let $(\omega^k)_{k\in \N}$ be the harmonic measures of their complementaries. Assume that $\ds \lim_k d(\Psi,\Psi^k)=0$.

We define for any $k$ a metric $d_k$ such that geometric cylinders for $\Psi^k$ and symbolic cylinders for $d_k$ have the same diameter.

We know that $d_k$ and $\omega^k$ fulfill the same properties as $d$ and $\omega$, with the same constants, because of our hypothesis : Annulus Condition and Bounded Contraction, and the fact that $(\Psi^k)$ is converging towards $\Psi$.

Let $s_k(X)=\diame{X}{d_k}$. Since $\ds\lim_k d(\Psi,\Psi^k)=0$, by Proposition \ref{contdiameters} we know that  $\ds\lim_k \mathcal{D}(s,s_k)=0$, and by Theorem \ref{compare} that $\ds\lim_k\mathcal{D}(\omega,\omega^k)=0$ we may thus apply Theorem \ref{contdim} to conclude for the continuity of the dimensions of harmonic measure.

\bibliographystyle{alpha}
\bibliography{biblio}

\end{document}